\documentclass[fleqn,12pt,twoside]{article}
\usepackage[headings]{espcrc1}
\usepackage{amsmath,amsthm,amssymb,slashbox,verbatim,graphicx}
\PassOptionsToPackage{ctagsplt,righttag}{amsmath}

\newcommand{\bN}{\mathbf{N}}

\newcommand{\bn}{\mathbf{n}}

\newcommand{\bx}{\mathbf{x}}

\newcommand{\cA}{\mathcal{A}}

\newcommand{\cL}{\mathcal{L}}

\newcommand{\tr}{\operatorname{tr}}
\newcommand{\grad}{\operatorname{grad}}

\newcommand{\beq}{\begin{equation}}
\newcommand{\eeq}{\end{equation}}
\newcommand{\bea}{\begin{eqnarray}}
\newcommand{\eea}{\end{eqnarray}}

\def\eg{\emph{e.g., }}
\def\ie{\emph{i.e., }}
\def\etal{\emph{et al.}}

\def\volfrac{\ensuremath{f}}

\newcounter{hours}
\newcounter{minutes}
\newcommand{\printtime}{\setcounter{hours}{\time/60}%
                        \setcounter{minutes}{\time-\value{hours}*60}%
\ifthenelse{\value{hours}<10}{0}{}\thehours:%
\ifthenelse{\value{minutes}<10}{0}{}\theminutes}
\advance\textheight by \topskip
\begin{document}
\bibliographystyle{unsrt}
%
%
%
%

\title{A Variational Level Set Approach for Surface Area Minimization 
of Triply Periodic Surfaces}

\author{Y. Jung\address[PRISM]{Princeton Institute for the Science and 
        Technology of Materials, Princeton University, Princeton, New Jersey 
        08544}\address[CHEM]{Department of Chemistry, Princeton University, 
        Princeton, New Jersey 08544},
        K. T. Chu\address[MAE]{Department of Mechanical and Aerospace 
        Engineering, Princeton University, Princeton, New Jersey 08544},
        and S. Torquato\addressmark[PRISM]\addressmark[CHEM]\address[PACM]
        {Program in Applied and Computational Mathematics, Princeton 
        University, Princeton, New Jersey 08544}\address[PCTP]{Princeton 
        Center for Theoretical Physics, Princeton University, Princeton, 
        New Jersey 08544}
}

\runtitle{Level Set Approach for Surface Area Minimization of Triply 
Periodic Surfaces}

\maketitle

\noindent \rule{6in}{1pt}

\begin{abstract}
In this paper, we study triply periodic surfaces with minimal surface area 
under a constraint in the volume fraction of the regions (phases) that the 
surface separates.  Using a variational level set method formulation, we 
present a theoretical characterization of and a numerical algorithm for 
computing these surfaces.  We use our theoretical and computational 
formulation to study the optimality of the Schwartz P, Schwartz D, and 
Schoen G surfaces when the volume fractions of the two phases are equal 
and explore the properties of optimal structures when the volume fractions 
of the two phases not equal.  Due to the computational cost of the fully, 
three-dimensional shape optimization problem, we implement our numerical 
simulations using a parallel level set method software package. 
\end{abstract}

\noindent \rule{6in}{1pt}

\section{Introduction}
\label{section1}
\begin{figure}[htb!]
\centerline{
\includegraphics[width=2.0in]{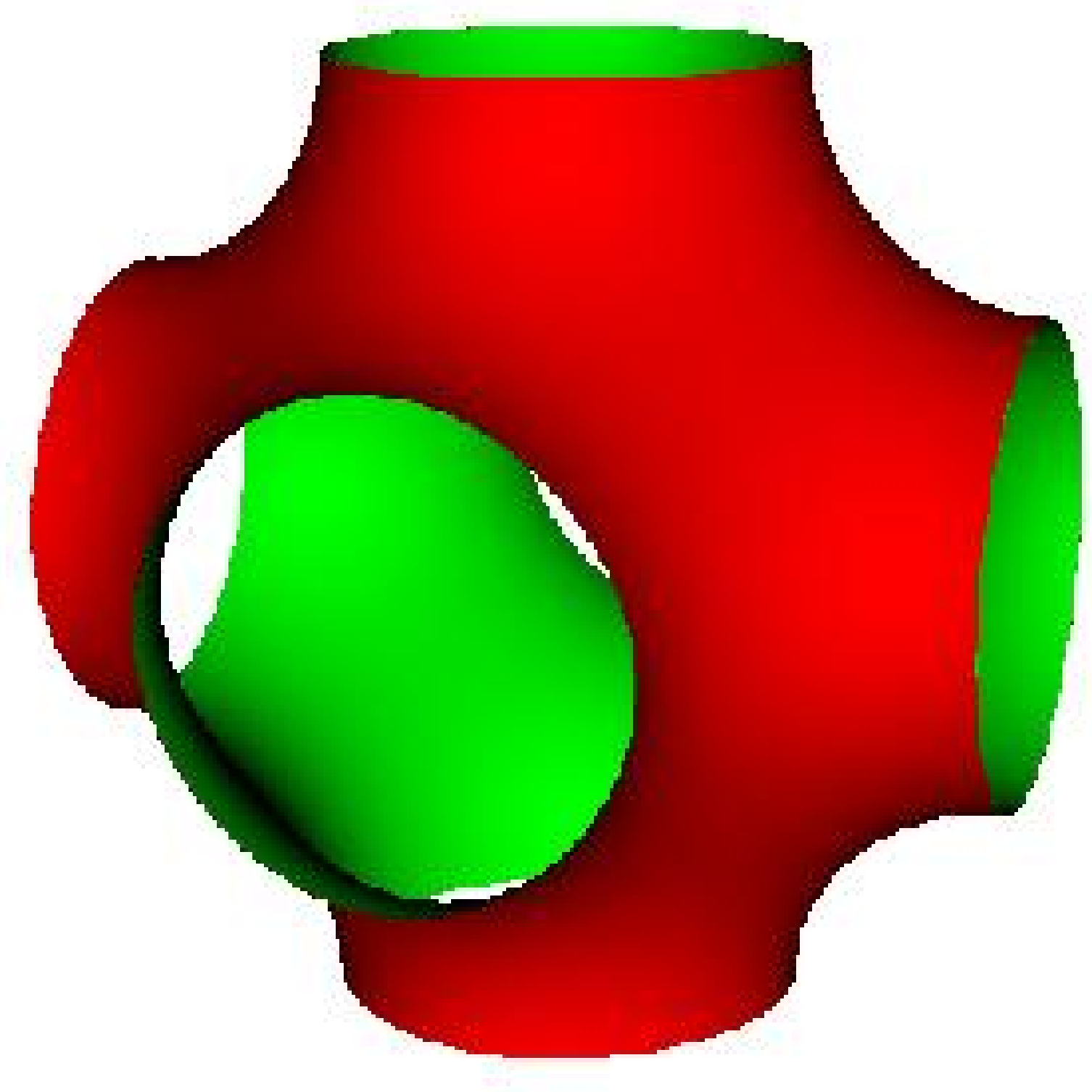} \hspace{0.0in}
\includegraphics[width=2.0in]{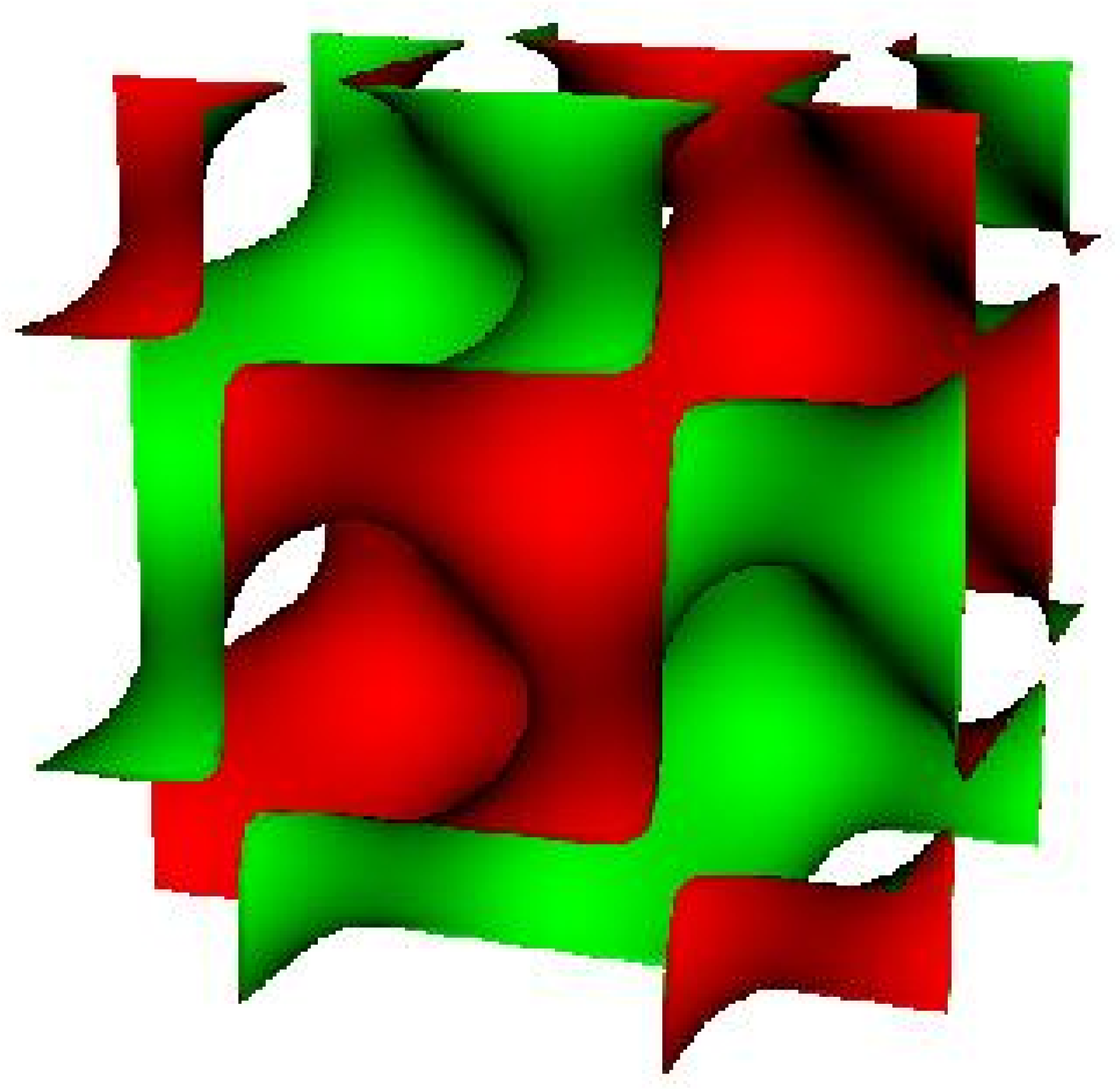} \hspace{0.0in}
\includegraphics[width=2.0in]{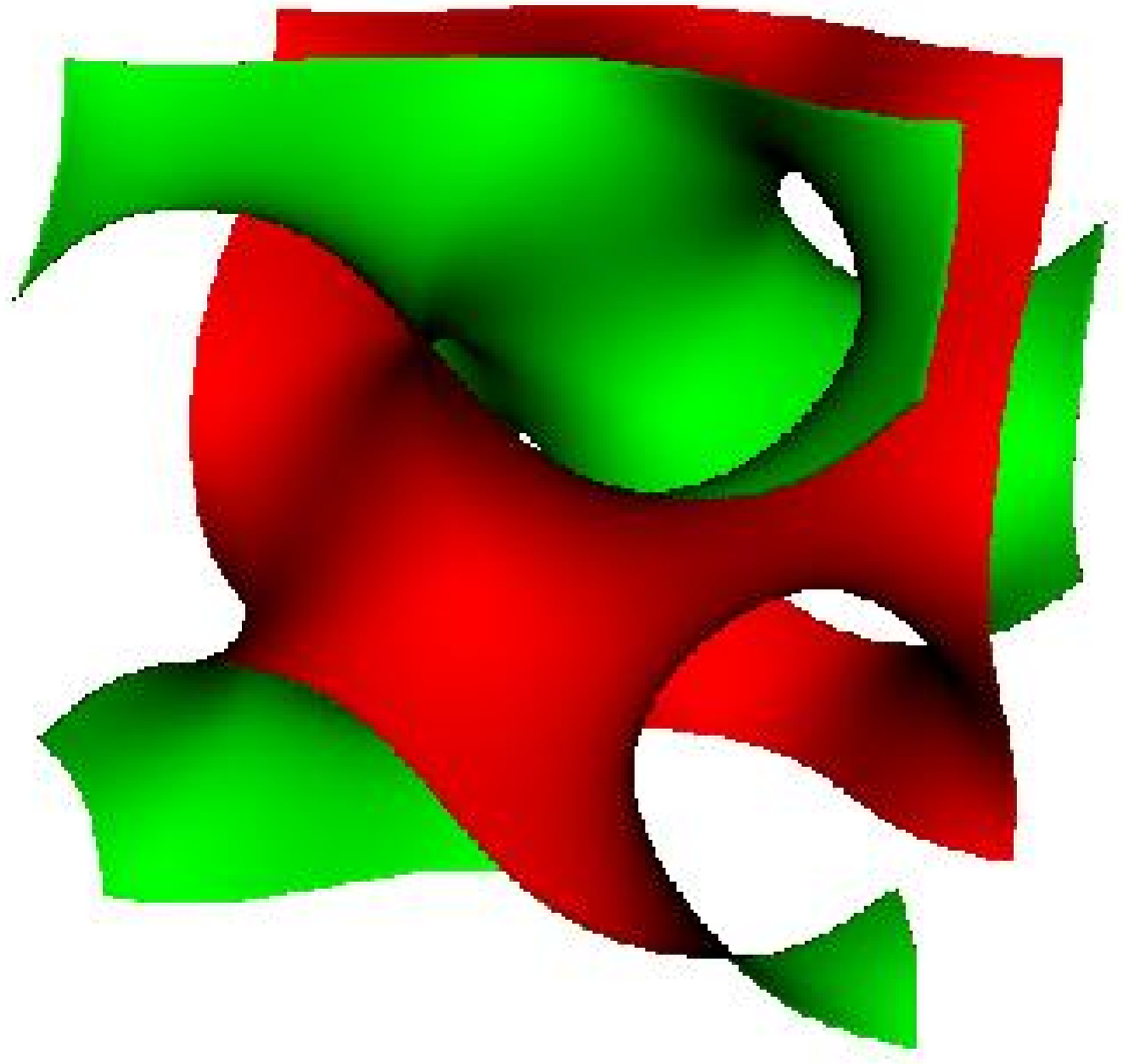}
}
\caption{Unit cell for three common minimal surfaces:
Schwartz P surface (left), Schwartz D surface (middle),
and Schoen G surface (right).}
\label{fig:tpms}
\end{figure}

Triply periodic minimal surfaces~\cite{Lord03,Mack93,Klin96} are objects of 
great interest to 
physical scientists, biologists, and mathematicians because they naturally 
arise in a variety of systems, including block copolymers~\cite{Ol98}, 
nanocomposites~\cite{Br01}, micellar materials~\cite{Zi00}, and lipid-water 
systems and certain cell membranes~\cite{Klin96,Gel94,Land95,Na96}.
There are two key feature of these surfaces:  (1) the mean 
curvature\footnote{The mean curvature, $H(\bx)$, at a point $\bx$ of a 
surface is the average of the two principal normal curvatures, $\kappa_1(\bx)$
and $\kappa_2(\bx)$: 
$H(\bx) = \frac{1}{2}\left(\kappa_1(\bx) + \kappa_2(\bx) \right)$.  $H(\bx)$
is also conveniently represented in terms of the divergence of the unit normal
vector (or equivalently, the trace of the gradient of $\bn$):
$H(\bx) = -\frac{1}{2} \nabla \cdot \bn = -\frac{1}{2} \tr( \grad \bn )$.}
is zero everywhere on the surface and (2) they are periodic in all three 
coordinate directions (i.e. they extend infinitely in all directions and 
possess the symmetry of one of the crystallographic space groups).  
An important subclass of triply periodic minimal surfaces are those that 
partition space into two disjoint but intertwining regions that are 
simultaneously continuous~\cite{An90}.  Examples of such surfaces
include the Schwartz primitive (P), the Schwartz diamond (D), and the 
Schoen gyroid (G) surfaces (see Figure~\ref{fig:tpms}).

Recently, there has been a resurgence of interest in triply periodic minimal 
surfaces because two-phase composites whose phases are separated by
such surfaces have been shown to be optimal with respect to multifunctional 
optimization of their material properties~\cite{To02d,To03a,To04a}.  
For instance, it has been shown that maximal values for the sum of the 
effective thermal and electrical conductivities is achieved by two-phase 
structures separated by Schwartz P and Schwartz D surfaces~\cite{To02d,To03a}.
These structures have also been discovered to be optimal for multifunctional
bulk modulus and electrical conductivity optimizations~\cite{To04a}.
As another example, the porous medium with the Schwartz P interface was 
found (via computer simulation) to have the largest 
fluid permeability over a range of microstructures examined~\cite{Jung05}.  

In~\cite{Jung05}, Jung and Torquato observed that the fluid 
permeabilities of porous, bicontinuous microstructures is inversely 
related to the \emph{total} interfacial surface area per unit volume.  This 
observation led them to conjecture that the maximal fluid permeability 
for a triply periodic porous medium is achieved by the structure that 
minimizes the total interfacial area~\cite{Jung05}.  One of the most 
interesting aspects of this conjecture is its focus on a \emph{global} 
property of the interface (\ie total surface area) rather than on 
\emph{local} (\eg differential) properties of the surface (\eg mean curvature).

Unfortunately, exploration of global properties of surfaces, such as the 
total surface area, seems to have received little attention in the 
literature.  Classically, the study of surfaces and their properties 
has been the realm of differential geometry.  In that field, research has 
traditionally focused on characterization and examination of minimal surfaces
(\ie surfaces with zero mean curvature).  The search for minimal surfaces 
has been ongoing activity since the mathematician H.~A. Schwartz published the 
first example of a minimal surface with full three-dimensional periodicity 
in 1865~\cite{schwarz65}. 
A systematic exploration of triply periodic surfaces with \emph{nonzero} mean 
curvature appears to have only been carried out relatively recently by 
Anderson, Davis, Scriven and Nitsche~\cite{An90}.  In their work, Anderson 
\etal~numerically computed surfaces of prescribed mean curvature (possibly 
non-constant) and studied their properties.  

In this paper, we consider triply periodic surfaces of nonzero mean curvature 
from a slightly different, but complementary, perspective.  Motivated by the 
fluid permeability conjecture and the scarcity of mathematical results on 
surfaces of nonzero mean curvature, we examine triply periodic surfaces that 
minimize the total interfacial surface area while satisfying a constraint on 
the volume fractions of the regions that the surface separates.  
Interestingly, we find that the optimal surfaces for this problem are
precisely those possessing a constant mean curvature.  This result allows
us to reproduce many of the results obtained by Andersen \etal.  However,
we emphasize that our approach comes from a fundamentally different 
perspective.  Whereas Andersen \etal~prescribe the mean curvature and consider
properties of the resulting surface, we specify the volume fraction and 
ask what the mean curvature of the surface must be to minimize the total 
surface area.
In general terms, we are interested in the properties of structures that 
arise when global geometric properties of a surface are prescribed (as 
opposed to examining the consequences of local properties).

The main conceptual foundations for our work come from the fundamental
ideas underlying the level set method~\cite{osher_book,seth99}.  Using 
these ideas, we have developed a novel numerical method for computing 
surfaces for which the surface area is a local minimum and that satisfy a 
specified volume fraction constraint.  The ideas underlying the level method 
were also used to provide a theoretical characterization of the surfaces 
that solve of our constrained optimization problem.
For both computational and theoretical purposes, the level set method 
formulation is convenient because it does not require explicit representation 
of the surface.  From a theoretical perspective, the level set method also 
has the advantage that the formulas for key geometric quantities, such as 
the surface area and volume, have a relatively simple form.

In section~\ref{section2}, we begin by developing a level set 
formulation for describing and analyzing triply periodic surfaces.  
In section~\ref{section3}, we use this framework to theoretically characterize 
extrema of the total interfacial area under a constraint on the volume fraction
of the phases.  In section~\ref{section4}, we present a numerical 
optimization procedure for computing surfaces that are local minima of the 
total interfacial area subject to a volume constraint. 
Finally, in section~\ref{section5}, we use a parallel implementation of our 
numerical algorithm to study several physically and mathematically interesting 
surfaces.  Parallelism was utilized to help deal with the computational cost 
of our fully three-dimensional optimization problems.

\section{Level Set Formulation}
\label{section2}
Our exploration of triply-periodic surfaces that are local minima of the 
total interfacial area subject to a volume constraint is based on ideas
underlying the application of the level set method to shape 
optimization~\cite{osher01,alexandrov_2005}.
Following the usual methodology~\cite{osher01,alexandrov_2005}, we represent 
surfaces as the zero level set of an embedding function, $\phi(x)$, defined 
throughout the three-dimensional volume and use variational calculus to 
derive the relationship between $\phi$ and global geometric properties of 
the surface.  

In this section, we present a detailed discussion of the level set formulation 
for surface area minimization in the presence of a volume fraction constraint.  
While it is clear that similar formulations have been used in previous 
studies~\cite{alexandrov_2005}, the explicit formulas and derivation for the 
most important geometric quantities were not provided.  In the present 
discussion, which is intended to fill this apparent gap in the literature, we 
derive formulas for the area and volume of a triply-periodic surface and the 
variations of these quantities with respect to $\phi$.  It is worth mentioning 
that the following derivation also holds for non-periodic surfaces and is 
easily extended to codimension-one surfaces in $n$-dimensional Euclidean space. 

\begin{figure}[htb!]
\centerline{
\includegraphics[width=2.1in]{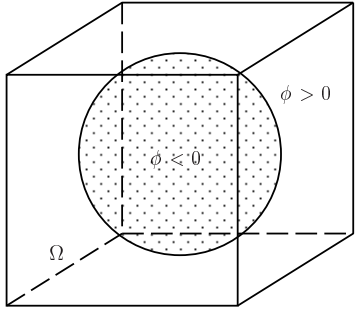}
}
\caption{Schematic diagram of a codimension-one surface and the sign of the 
level set function, $\phi$, in various regions of the domain.  The shaded 
and unshaded regions correspond to phase 1 and phase 2, respectively.}
\label{fig:schematic_diagram}
\end{figure}

\begin{figure}[htb!]
\centerline{
\includegraphics[width=2.0in]{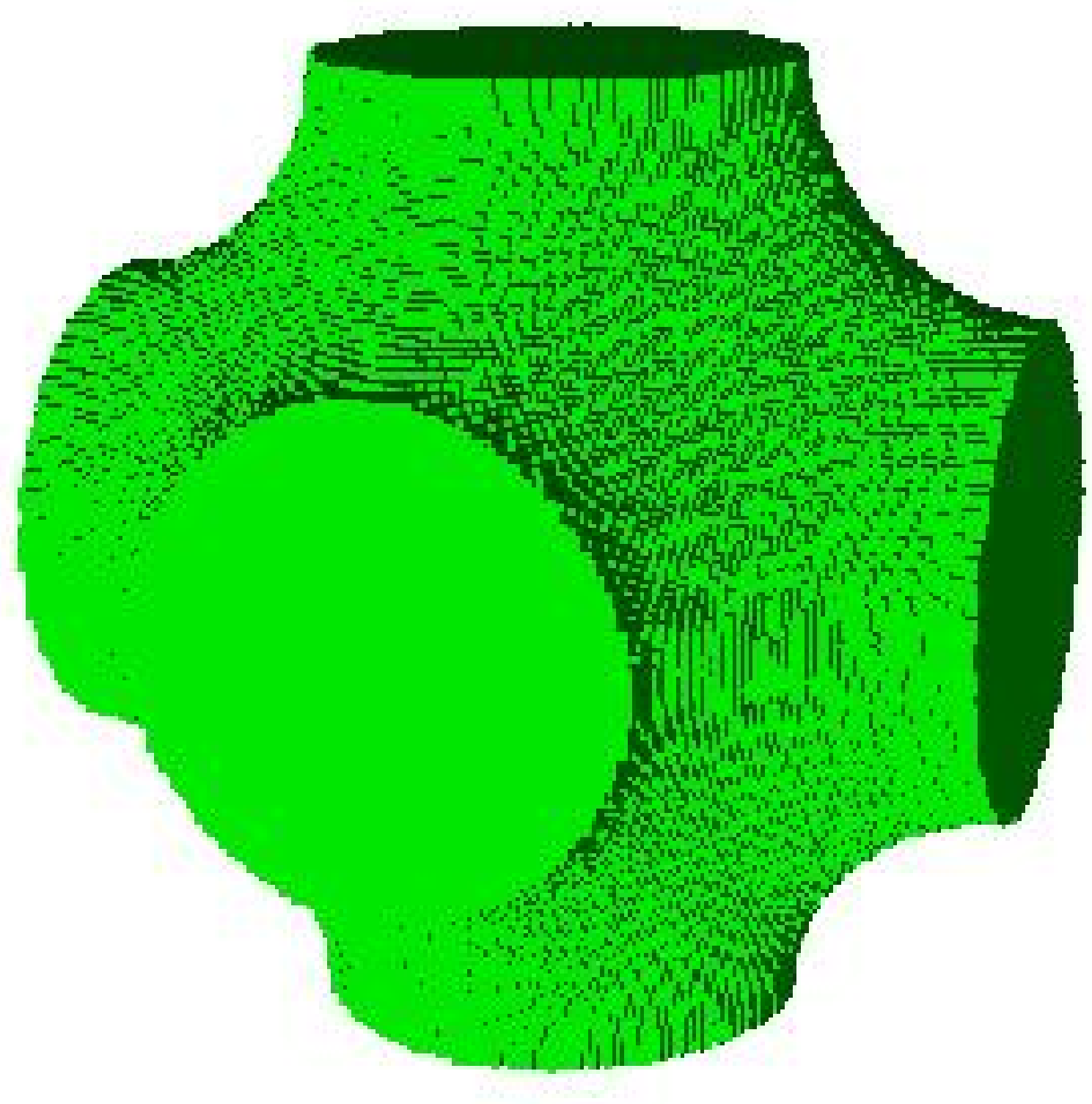} \hspace{0.0in}
\includegraphics[width=2.0in]{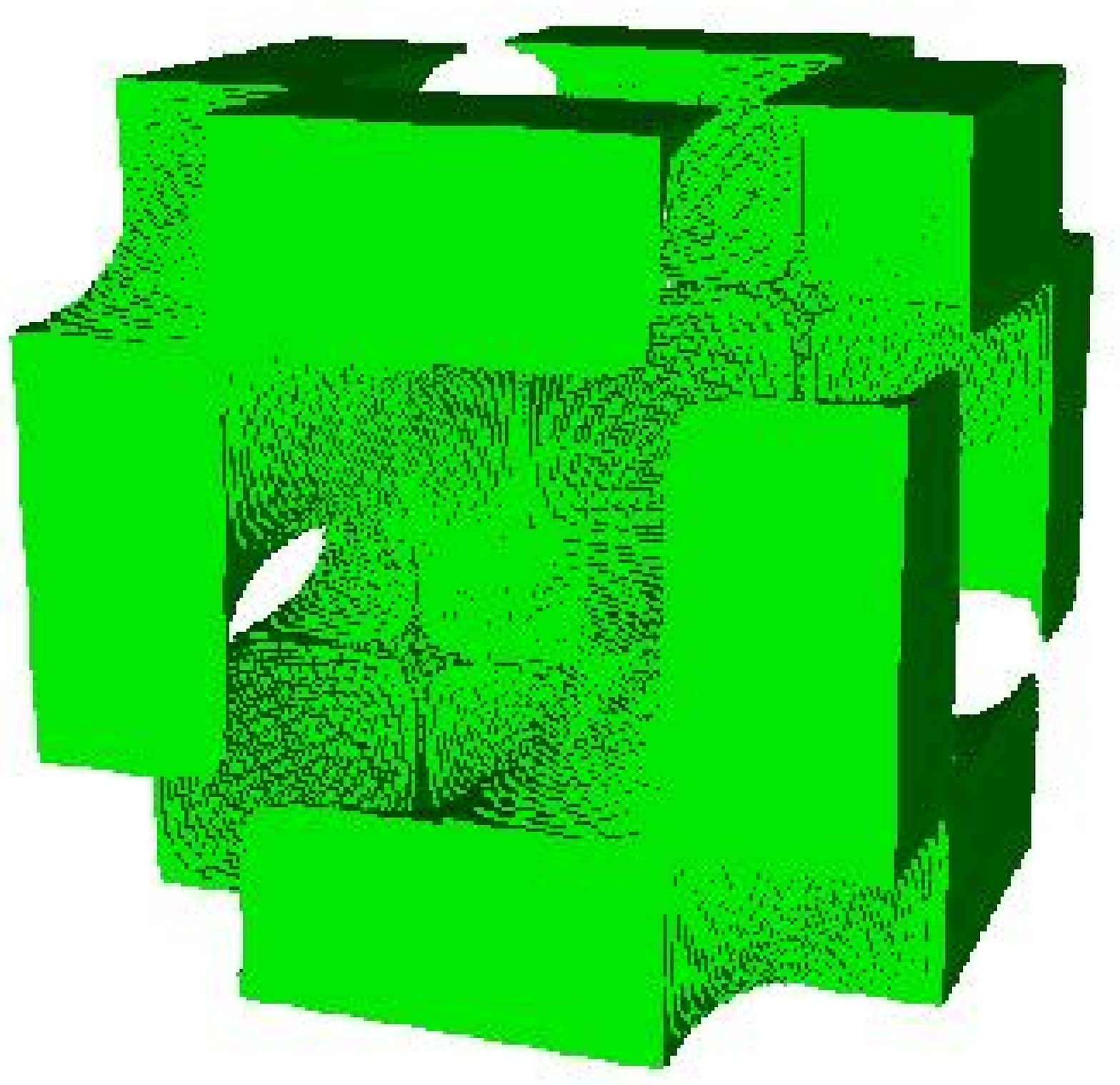} \hspace{0.0in}
\includegraphics[width=2.0in]{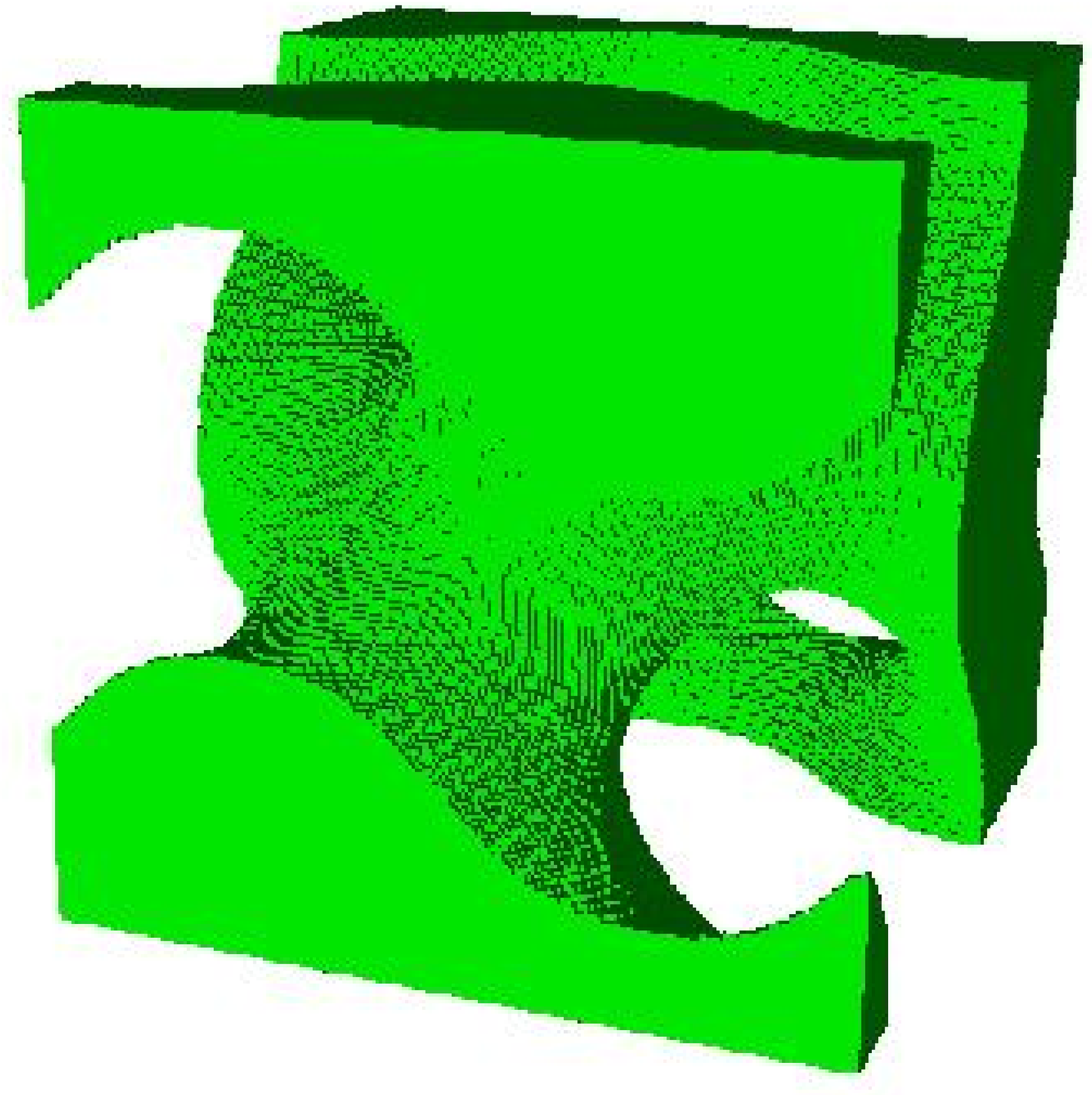} 
}
\caption{The phase 1 region corresponding to three different choices for 
$\Gamma$:  Schwartz P surface (left), Schwartz D surface (middle),
and Schoen G surface (right).}
\label{fig:phase_periodicity}
\end{figure}

For the triply periodic problems that we are interested in, the domain is 
taken to be a unit cell of the periodically repeating structure.  Let the 
surface of interest, $\Gamma$, be represented by the zero level set of $\phi$.  
Then $\Gamma$ divides the unit cell into two distinct phases.  Without loss 
of generality, we define the region where $\phi(\bx) < 0$ to be phase $1$ 
(see Figures~\ref{fig:schematic_diagram} and~\ref{fig:phase_periodicity}).
An important assumption in our work is that a shift of the unit cell by a 
lattice vector does \emph{not} cause an interchange of the phases.  In this 
paper, systems possessing this property will be called \emph{phase-periodic}.

In terms of the embedding function, $\phi$, our minimization problem may be 
stated as
\begin{equation}
\textrm{minimize}\ \cA(\phi)\ \textrm{subject to}\ \volfrac(\phi) = \volfrac_o,
\label{eq:optimization_problem}
\end{equation}
where $\cA(\phi)$ is the total surface area of the zero level set, 
$\volfrac(\phi)$ is the volume fraction of phase $1$, and $\volfrac_o$ is 
the desired volume fraction for phase $1$.  Unfortunately, there are 
currently no theoretical results concerning the global minimum of our 
problem.  However, as we shall see, there are many local minima for this 
problem.  Because many interesting surfaces arise as local minima, we shall 
focus our attention on understanding the structure of these surfaces 
and neglect, for the moment, the question of global minimality. 

We use the method of Lagrange multipliers to solve the optimization problem 
(\ref{eq:optimization_problem}). The Lagrangian is given by
\begin{equation} 
  \cL(\phi,\lambda) = 
  \cA(\phi) + \lambda\ \left( \volfrac(\phi)-\volfrac_o \right),
\end{equation} 
where $\lambda$ is the Lagrange multiplier.
Taking the variation of $\cL$ with respect to $\phi$, we find that a 
necessary condition for a minimizer is 
\begin{equation}
  \delta \cL(\phi,\lambda) = 
  \delta \cA(\phi) + \lambda\ \delta \volfrac(\phi)\ =\ 0.
\end{equation}
This condition, together with the volume fraction constraint
\begin{equation} 
  \volfrac(\phi) - \volfrac_o= 0
\end{equation} 
allows us to compute $\phi$ and $\lambda$. 

\subsection{Volume integral formulation of surface integrals}
\label{subsection2.1}
To compute the variations of $\cA$ and $\volfrac$ required by our Lagrangian
formulation of the optimization problem, it is convenient to first derive 
a relation that allows us to convert between surface integrals over $\Gamma$
and volume integrals over entire unit cell.  
Let $\Omega$ denote the entire unit cell, and let $\bN$ denote an outward
pointing unit normal vector on the boundary of the unit cell, 
$\partial \Omega$.
Similarly, let $\Omega_1$ and $\bn$ denote the region corresponding to 
phase $1$ and an outward pointing normal vector on 
$\Gamma = \partial \Omega_1$.
Now, consider the surface integral of an arbitrary function $p(\phi)$ over
$\Gamma$:
\begin{eqnarray}
\int_{\Gamma} p(\phi)\ dS.
\end{eqnarray}
Applying the divergence theorem, we can rewrite this surface integral 
as a volume integral over phase $1$ and surface integrals over the 
intersections of phase $1$ with the boundaries of the unit cell:
\begin{eqnarray}
\int_{\Gamma} p(\phi)\ dS &=& \int_{\Gamma} p(\phi) \bn \cdot \bn\ dS 
          = \int_{\Gamma} 
   \left( p(\phi) \frac{\nabla \phi}{\left\| \nabla \phi \right\|} \right) 
   \cdot \bn\ dS \nonumber \\
          &=& \int_{\Omega_1} 
   \nabla \cdot \left( p(\phi) \frac{\nabla \phi}{\left\| \nabla \phi \right\|} 
   \right)\ dV - \int_{B} 
   \left( p(\phi) \frac{\nabla \phi}{\left\| \nabla \phi \right\|} \right) 
   \cdot \bN\ dS, \label{eq:integral_conversion_intermediate}
\end{eqnarray}
where $B$ is the union of the intersections of phase $1$ with the boundaries
of the unit cell and we have identified $\bn$ with 
$\frac{\nabla \phi}{\left\| \nabla \phi \right\|}$. 
Notice that the last surface integral term vanishes because we have assumed 
that $\phi$ is phase-periodic which means that
$\left( p(\phi) \frac{\nabla \phi}{\left\| \nabla \phi \right\|} \right)$ 
is invariant under a translation of $\bx$ by any lattice vector.
Using phase-periodicity, we can cancel out contributions of surface integrals 
from opposite faces of the unit cell because $\bN$ on opposite faces have
opposite sign.  

Next, we observe that the volume integral over phase $1$ can be rewritten as 
a volume integral over the entire unit cell:
\begin{equation} 
\int_{\Omega_1} \nabla \cdot \left( p(\phi) \frac{\nabla \phi}
{\left\| \nabla \phi \right\|} \right)\ dV
 = \int_\Omega
  \left ( 1 - \Theta(\phi) \right)
   \nabla \cdot \left( p(\phi) \frac{\nabla \phi}{\left\| \nabla \phi \right\|} 
   \right) dV,
\end{equation} 
where $\Theta(\phi)$ is the Heaviside function  
$$
\Theta(\phi)\ =\
 \begin{cases}
   1 & \qquad \phi\ \ge\ 0 \\
   0 & \qquad \phi\ <\ 0. 
 \end{cases}
$$
This expression can be further simplified via an integration by parts
procedure to obtain:
\begin{eqnarray}
\int_{\Gamma} p(\phi)\ dS 
        &=& -\int_\Omega \nabla \left ( 1 - \Theta(\phi) \right) \cdot
            \left( 
               p(\phi) \frac{\nabla \phi}{\left\| \nabla \phi \right\|} 
            \right ) dV \nonumber \\
        &\ & +\ \int_{\partial\Omega}
            \left ( 1 - \Theta(\phi) \right) 
            \left( 
               p(\phi) \frac{\nabla \phi}{\left\| \nabla \phi \right\|} 
            \right )
            \cdot \bN\ dS \nonumber \\
        &=& \int_\Omega \delta(\phi) \nabla \phi \cdot
            \left( 
               p(\phi) \frac{\nabla \phi}{\left\| \nabla \phi \right\|} 
            \right ) dV \nonumber \\
        &=& \int_\Omega p(\phi)\ \delta(\phi) \left\| \nabla \phi \right\| dV.
            \label{eq:integral_conversion}
\end{eqnarray}
where $\delta(\phi)$ is the Dirac delta function and we have eliminated
the surface integral contribution by appealing to a similar line of reasoning 
as used when eliminating the surface integral term from
(\ref{eq:integral_conversion_intermediate}).
The form of (\ref{eq:integral_conversion}) is not unexpected -- the
delta function ensures that only values of $p(\phi)$ on the surface $\Gamma$ 
contribute to the integral and the $||\nabla \phi||$ is the Jacobian
that arises when $\phi$ (as opposed to a spatial coordinate in the direction
normal to the surface) is used as the argument of the delta function.
Equation (\ref{eq:integral_conversion}) forms the foundation for much of 
the following mathematical development.  

\subsection{Variation of the surface area}
\label{subsection2.2}
Using the formula for the area of $\Gamma$, $\cA(\phi) = \int_\Gamma dS$, and 
(\ref{eq:integral_conversion}), we find that
\beq
\label{eq:surface}
  \cA(\phi) = \int_\Omega \delta(\phi) \left\| \nabla \phi \right\|\ dV 
\eeq
Taking the variation of this expression with respect to $\phi$ yields 
\bea
\delta \cA(\phi) &=& \int_\Omega 
  \left [ 
    \delta \left\{ \delta (\phi) \right\} \left\| \nabla \phi \right\|
  + \delta (\phi) \delta \left\| \nabla \phi \right\|
  \right ] dV  \nonumber \\
         &=& \int_\Omega 
  \left [ 
    \delta' (\phi) \left\| \nabla \phi \right\| \delta \phi
  + \delta (\phi) \frac{\nabla \phi \cdot \nabla (\delta \phi) }
                       {\left\| \nabla \phi \right\|}
  \right ] dV  \nonumber \\
         &=& \int_\Omega 
  \left [ 
    \nabla (\delta (\phi)) \cdot \frac{\nabla \phi}{\left\| \nabla \phi \right\|}
       \delta \phi
  + \delta (\phi) \frac{\nabla \phi \cdot \nabla (\delta \phi) }
                       {\left\| \nabla \phi \right\|}
  \right ] dV 
\label{eq:delta_F_intermediate}
\eea
where $\delta^\prime(\phi)$ is the first derivative of the Dirac delta function
and we have used the relationship $\nabla \delta (\phi) = \delta^\prime(\phi) 
\nabla \phi$. We can further simplify $\delta \cA(\phi)$ by using the product 
rule to make the following substitutions:
\beq
  \nabla (\delta (\phi)) \cdot \frac{\nabla \phi}{\left\| \nabla \phi \right\|}
  = \nabla \cdot \left( \delta(\phi) 
       \frac{\nabla \phi}{\left\| \nabla \phi \right\|} \right)
  - \delta(\phi) \nabla \cdot \left( 
      \frac{\nabla \phi}{\left\| \nabla \phi \right\|} \right)
\eeq
and
\beq
  \delta (\phi) \frac{\nabla \phi \cdot \nabla (\delta \phi) }
                       {\left\| \nabla \phi \right\|} 
  = \nabla \cdot \left( \delta \phi \delta(\phi) 
      \frac{\nabla \phi}{\left\| \nabla \phi \right\|} \right)
  - \delta \phi \nabla \cdot \left( \delta(\phi) 
      \frac{\nabla \phi}{\left\| \nabla \phi \right\|} \right).
\eeq
Thus, (\ref{eq:delta_F_intermediate}) can be rewritten as
\bea
\delta \cA(\phi) &=& \int_\Omega 
  \left[ \nabla \cdot \left( \delta \phi \delta(\phi) 
      \frac{\nabla \phi}{\left\| \nabla \phi \right\|} \right)
  - \delta \phi \delta(\phi) \nabla \cdot \left( 
      \frac{\nabla \phi}{\left\| \nabla \phi \right\|} \right)
  \right] dV  \nonumber \\
  &=& - \int_\Omega \delta \phi \delta(\phi) \nabla \cdot \left( 
        \frac{\nabla \phi}{\left\| \nabla \phi \right\|} \right) dV,
\label{eq:delta_F_final_volume}
\eea
where we have again used phase-periodicity to eliminate the surface integral 
over the unit cell that arises when applying the divergence theorem to the
first integrand.  
Finally, we may convert $\delta \cA$ back into a surface integral over 
$\Gamma$ by comparing (\ref{eq:delta_F_final_volume}) with 
(\ref{eq:integral_conversion}): 
\begin{equation}
\delta \cA(\phi) = -\int_{\Gamma} (\nabla \cdot \bn)\ \frac{\delta \phi}
{\|\nabla \phi\|}\ dS.
\end{equation}

\subsection{Variation of the volume fraction constraint}
\label{subsubsection2.3}
To compute the variation of the volume fraction constraint, we begin by 
writing the definition of the volume fraction as an integral over the 
entire unit cell: 
\begin{equation}
\volfrac(\phi) = \int_{\Omega_1} 1\ dV \label{eq:volume_frac_formula} 
        = \int_\Omega \left ( 1 - \Theta(\phi) \right) dV.
\end{equation}
Note that in writing (\ref{eq:volume_frac_formula}), we have implicitly 
assumed that the volume of the unit cell is $1$.
Next, we compute the variation of $\volfrac(\phi)$ as follows: 
\begin{equation}
\delta \volfrac(\phi) = \int_{\Omega} \delta \left ( 1 - \Theta(\phi) \right)
  dV \\
  = -\int_\Omega \delta(\phi) \delta \phi dV \\ 
  = -\int_{\Gamma} \frac{\delta \phi}{\|\nabla \phi\|} dS,
\label{eq:delta_G}
\end{equation}
where we have again made use of (\ref{eq:integral_conversion}) to convert 
the volume integral over the entire unit cell to a surface integral over 
$\Gamma$.

\section{Theoretical Analysis of Minimization Problem}
\label{section3}

\subsection{Characterization of local extrema}
\label{subsection3.1}
Using the machinery developed in the previous section, we can characterize 
the local extrema of the optimization problem (\ref{eq:optimization_problem}).
A straightforward application of the calculus of variations with constraints 
leads to the following
\newtheorem*{extrema_total_surface_area}
  {Theorem\label{thm:extrema_total_surface_area}}
\begin{extrema_total_surface_area}
$\cA(\phi)$ is a local extremum subject to the constraint 
$\volfrac(\phi) = \volfrac_o$ 
if and only if the surface $\Gamma$ has constant mean curvature.
\end{extrema_total_surface_area}

\begin{proof}
$(\Rightarrow)$  Suppose that $\Gamma$ is a local extremum of the 
total surface area subject to the specified volume fraction constraint.  
Then it is a local extremum (without constraint) of the Lagrangian 
\cite{hildebrand_book}:
\beq
\cL(\phi, \lambda) \equiv \cA(\phi) + \lambda 
\left( \volfrac(\phi) - \volfrac_o \right),
\eeq
where $\lambda$ is a scalar Lagrange multiplier.
Using the expressions for the variation of the surface area and volume 
fraction derived in the previous section, the total variation of 
$\cL(\phi,\lambda)$ is given by
\bea
\delta \cL(\phi, \lambda) &=& \delta \cA(\phi) 
  + \delta \left [ \lambda \left ( \volfrac(\phi) - \volfrac_o \right ) 
           \right ] \nonumber \\
  &=& -\int_{\Gamma} \left [ (\nabla \cdot \bn) + \lambda \right ]
                         \frac{\delta \phi}{\| \nabla \phi \|} dS
  + \left ( \volfrac(\phi) - \volfrac_o \right ) \delta \lambda.
  \label{eq:variation_L}
\eea
Since this expression must be of zero for \emph{any} variation
in $\phi$ and $\lambda$, (\ref{eq:variation_L}) implies
that
\bea
  \nabla \cdot \bn + \lambda &=& 0 \ \ \ \textrm{on\ } \Gamma
  \label{eq:curvature_constraint} \\
  \volfrac(\phi) - \volfrac_o &=& 0.
  \label{eq:volume_constraint}
\eea
Equation (\ref{eq:volume_constraint}) is just the volume fraction constraint
(which is expected and arises for all applications of the method of 
Lagrange multipliers).  Equation (\ref{eq:curvature_constraint}),
however, leads to the conclusion that the mean curvature 
is constant over the entire surface and is equal to half 
the Lagrange multiplier.

$(\Leftarrow)$  Now, suppose that $\Gamma$ is a constant mean
curvature surface.  Then, the variation in the area simplifies to yield
\bea 
  \delta \cA(\phi) 
  = -\int_{\Gamma} (\nabla \cdot \bn) \frac{\delta \phi}
     {\| \nabla \phi \|} dS
  = - (\nabla \cdot \bn) \int_{\Gamma} \frac{\delta \phi}
     {\| \nabla \phi \|} dS
  = (\nabla \cdot \bn) \delta \volfrac(\phi).
\eea 
Since any variation of the level set function $\phi$ that does not change
the volume fraction must satisfy $\delta \volfrac(\phi) = 0$, we see that 
$\delta \cA(\phi) = 0$ for any volume fraction preserving variation in $\phi$.  
In other words, $\Gamma$ is a local extremum of the total surface area 
subject to a volume fraction constraint. 
\end{proof}

Before moving on, it is worth mentioning some special cases of local 
extrema of the total surface area with a constrained volume fraction. 
First, there are the minimal surfaces.  Because, they have zero mean 
curvature, the theorem immediately implies that all minimal surfaces are 
local extrema of the total surface area under a volume fraction constraint.  
However, minimal surfaces also have the interesting property that they are 
local extrema of the total surface area \emph{without} any constraint.  
This result follows directly from the expression for the variation in the 
total surface area because the integrand is identically zero for minimal 
surfaces.
A second class of important examples is the class of constant, nonzero mean 
curvature surfaces.  Two examples are the sphere and infinite cylinder.
Both of these objects have a constant, nonzero mean curvature, and thus are 
local extrema of the total surface area when the volume is constrained.  
Further examples are provided by any of the nonzero mean curvature surfaces
studied by Anderson \etal~\cite{An90}.  

These examples demonstrate that the class of minimal surfaces and the class
of local extrema of the total surface area under a volume fraction constraint
are truly distinct.  That surfaces arising from an analysis of local 
geometric properties differ from those arising when global geometric 
properties are studied underscores the importance of considering the global 
geometric features of surfaces.  These considerations are especially
important from a physical perspective; while the driving force for surface 
evolution is often local in nature, global constraints (\eg conserved 
quantities) are almost always present and can affect the global structures 
that arise.

\section{Numerical Optimization Procedure}
\label{section4}
While Theorem~\ref{thm:extrema_total_surface_area} provides a theoretical
characterization of local minima of the total surface under a constrained
volume fraction, it is does not provide a means for obtaining optimal 
surfaces.  In this section, we present a numerical procedure for computing 
locally optimal surfaces.  Our approach follows the work of Osher and 
Santosa~\cite{osher01} which evolves the embedding function $\phi$ along 
steepest descent directions using the evolution equation:
\begin{equation}
\delta \phi + v(\bx) ||\nabla \phi|| = 0,
\label{eq:HamiltonJacobi}
\end{equation}
with $v(\bx)$ chosen so that $\phi$ remains in the set of embedding
functions that satisfy the volume fraction constraint.
As pointed out in~\cite{osher01}, this equation is equivalent to a 
Hamilton-Jacobi equation if the change in $\phi$ is viewed as occurring
continuously in time.  In addition, we use an auxiliary Newton iteration 
to explicitly enforce the volume fraction constraint when numerical error 
causes the volume fraction to drift beyond an acceptable tolerance and to
generate initial structures that satisfy a prescribed volume fraction 
constraint~\cite{osher01}.  

\subsection{Descent direction and projected gradient algorithm}
\label{subsection4.1}
Following~\cite{osher01}, we choose the velocity field $v(\bx)$ to be the 
steepest descent direction for the Lagrangian, $\cL$.  
From (\ref{eq:variation_L}), we know that the variation in $\cL$ 
is given by
\begin{equation}
\label{eq:descent}
\delta \cL(\phi,\lambda) = 
 -\int_{\Gamma} \left[ (\nabla \cdot \bn)\ + \lambda \right]
                \frac{\delta \phi}{\|\nabla \phi\|} dS.
\end{equation}
Therefore, choosing $\delta \phi$ such that
\begin{equation}
\label{eq:delta_phi}
 \delta \phi = \left[ (\nabla \cdot \bn)\ + \lambda \right] \|\nabla \phi\|
\end{equation}
ensures that $\cL$ decreases at every iteration.
Comparing with (\ref{eq:HamiltonJacobi}), we can identify the velocity 
field $v(\bx)$ as 
\begin{equation}
\label{eq:normal_velocity}
v(\bx) =  -\left[ (\nabla \cdot \bn)\ + \lambda \right].
\end{equation}
This velocity field has a natural interpretation: the motion of the interface
is driven by mean curvature (which tends to shrink the local area of the 
interface) and the Lagrange multiplier (which tries to keep the volume
fraction from changing).

To compute the value of $\lambda$ to use for each iteration, we insist that 
$\delta \phi$ is chosen to satisfy
\begin{equation}
  \delta \volfrac(\phi) = 0.
  \label{eq:numerical_volume_constraint}
\end{equation}
Essentially, we are ensuring that $\phi + \delta \phi$ satisfies the 
constraint equation linearized about the current iteration of $\phi$:
\beq
 \volfrac(\phi + \delta \phi) \approx \volfrac(\phi) + \delta \volfrac(\phi)
\eeq
Substituting (\ref{eq:delta_phi}) into (\ref{eq:delta_G}) and enforcing
(\ref{eq:numerical_volume_constraint}), we find that the Lagrange multiplier 
should be set equal to 
\begin{equation}
\label{eq:Lagrange_multiplier}
\lambda = -\frac{\int_{\Gamma} (\nabla \cdot \bn) dS}
                {\int_{\Gamma} 1\ dS} =
          -\frac{\int_{\Gamma} (\nabla \cdot \bn) dS}
                {\cA(\phi)}. 
\end{equation}
Notice that $\lambda$ is equal to twice the average value of the 
the mean curvature over the surface.

\subsubsection{Extension of the velocity off of the interface}
The velocity of points off of the zero level set are obtained by 
extending the velocity on $\Gamma$ in the normal 
direction~\cite{osher_book}.  We decided to set the velocity in this way
to avoid the mathematical singularities that arise for non-zero level sets 
with very large mean curvatures when using the ``natural'' 
velocity extension\footnote{The ``natural'' velocity extension uses
(\ref{eq:normal_velocity}) for points off of the zero level set with 
$\bn$ replaced by $\frac{\nabla \phi}{||\nabla \phi||}$.}.  
These singularities, which appear at cusps of the embedding function $\phi$,
lead to numerical difficulties (\eg strict stability constraint on the 
effective time step size).

\subsection{Newton iteration to enforce volume fraction constraint}
\label{subsection4.2}

Because we are only approximately enforcing the constraint at each time step,
the iterates of the embedding function $\phi$ will eventually fail to satisfy
the volume fraction constraint.  To put an iterate back onto the feasible 
set after the constraint has been violated by more than an acceptable
tolerance, we use a Newton iteration as in~\cite{osher01}.  
We also use the Newton iteration to ensure that the initial condition
for the embedding function satisfies the desired volume fraction constraint.
For situations where the target volume fraction constraint is
far from the actual volume fraction of the initial iterate, we use 
continuation in $\volfrac_o$ to improve the convergence of the volume 
fraction constraint algorithm~\cite{boyd_book}. 

We begin by considering the volume fraction of a corrected embedding function:
$\volfrac(\phi^{(0)} + \delta \phi)$, where $\phi^{(0)}$ is the uncorrected 
embedding function.  Next, we think of $\delta \phi$ in
(\ref{eq:delta_phi}) as a function of $\lambda$ (taken to be an unknown) and 
introduce a scale factor $\alpha > 0$:
\beq
\delta \phi(\lambda) = 
  \alpha \left[ (\nabla \cdot \bn)\ + \lambda \right] \|\nabla \phi\|.
\label{eq:delta_phi_newton}
\eeq
This choice for $\delta \phi(\lambda)$ was originally proposed by Osher and 
Santosa in~\cite{osher01}.  While there are certainly other possible choices
for $\delta \phi(\lambda)$, we opted to use (\ref{eq:delta_phi_newton}) 
because it tends to preserve the shape of the surface, it performed 
well in the present work, and it allowed us to reuse code written for 
other portions of the computation.

Note that (\ref{eq:delta_phi_newton}) allows us to think of 
$\volfrac(\phi^{(0)} + \delta \phi)$ as a function of $\lambda$.
We can then use a Newton iteration to compute the value of $\lambda$ 
(and therefore the appropriate correction $\delta \phi$) required so 
that the volume fraction constraint is satisfied: 
$\volfrac(\phi^{(0)} + \delta \phi) = 0$. 
It is important to recognize that the choice $\delta \phi$ in the Newton 
iteration is different from the choice of $\delta \phi$ in the optimization 
iteration (\ref{eq:delta_phi}).  In the latter case, $\lambda$ has a known 
value and can be explicitly computed; in the former case, $\lambda$ is an 
unknown variable that is determined through a Newton iteration.

Taking advantage of the level set formulation of our problem, calculating the
derivative of $\volfrac(\lambda)$ with respect to $\lambda$ is 
straightforward:
\begin{eqnarray}
D_\lambda \volfrac(\phi^{(0)} + \delta \phi(\lambda)) &=& 
  -\int_{\Omega} D_\lambda H(\phi^{(0)} + \delta \phi)  dV \nonumber \\
  &=& -\int_\Omega D_{\phi^{(0)}+\delta \phi} 
         H(\phi^{(0)} + \delta \phi) D_\lambda 
         \left(\phi^{(0)} + \delta \phi \right) dV \nonumber \\ 
  &=& -\alpha \int_\Omega \delta(\phi^{(0)} + \delta \phi) 
         \left\| \nabla \phi \right\| dV.
\label{eq:d_volfrac_d_lambda}
\end{eqnarray}

\subsubsection{Choice of scale factor}
The scale factor is necessary to ensure that the zero level set of
$\left( \phi^{(0)} + \delta \phi \right)$ does not vanish and remains within 
the unit cell; otherwise, the derivative of $\volfrac$ with respect $\lambda$ 
which is computed via as integral over the unit cell in 
(\ref{eq:d_volfrac_d_lambda}) may vanish.  
Intuitively, the scale factor helps make the discretized version of 
$\delta \phi$ a reasonable approximation to an infinitesimal variation
of $\phi$.  Without it, the discrete approximation to $\delta \phi$
would be equal to $\left[ (\nabla \cdot \bn) + \lambda \right]$, which
may not be small for all Newton iterations.
In our simulations, we employed a scaling factor of 
$\alpha = (\min \{\Delta x_1, \Delta x_2, \Delta x_3\})^2$ where 
$\Delta x_i$ is the grid spacing in each $i$-th coordinate direction.
This choice for $\alpha$ was selected by numerical experimentation.

\subsection{Numerical implementation issues}
\label{subsection4.3}
The optimization algorithm based on the algorithmic components discussed in 
the previous two sections may be summarized in the following steps:
\begin{enumerate}
\item Generate an initial configuration for $\phi_{0}(\bx)$. 
    
\item Evolve the level set function to time $t_{n+1}$ if $\Delta \cA_n >$ tol.

  \begin{enumerate}
  \item If necessary, use Newton iteration to enforce the volume fraction 
        constraint. 
     \begin{enumerate}
     \item Select an initial guess for $\lambda$.

       \item Repeat the following steps until convergence:
       \begin{enumerate}

         \item Compute $\delta \phi$ using (\ref{eq:delta_phi_newton}). 

         \item Compute $\phi^{(0)} + \delta \phi$.

         \item Update $\lambda$ using:
         $  \lambda^{(k+1)} = \lambda^{(k)} - \frac
        {\volfrac \left( \phi^{(0)}+\delta \phi(\lambda^{(k)}) \right) 
          - \volfrac_o}
        {D_\lambda \volfrac 
           \left( \phi^{(0)}+\delta \phi(\lambda^{(k)}) \right)}
         $
           where ($k$) is a $k$-th Newton iterate.
       \end{enumerate}

     \end{enumerate}

  \item Compute the Lagrange multiplier $\lambda_n$ using
        (\ref{eq:Lagrange_multiplier}).

  \item Compute the descent direction $\delta \phi_n$ using
        (\ref{eq:delta_phi}).

  \item Update $\phi$ using 
        $\phi_{n+1}(\bx)$ to $\phi_n(\bx)+ \beta \delta \phi_n(\bx)$.
        $\beta$ is a scale factor required for stability of the numerical
        scheme. 
        It is equivalent to a stable time step size when 
        (\ref{eq:HamiltonJacobi}) is viewed as a Hamilton-Jacobi equation 
        varying continuously in time.

  \item Periodically reinitialize $\phi(x)$ to an approximate distance 
        function within a sufficiently wide band around the zero level
        set of $\phi(x)$.

  \end{enumerate}

\end{enumerate}

\subsubsection{Synergy of reinitialization and 
               volume fraction constraint enforcement}
Both reinitialization and enforcement of the volume fraction constraint
introduce errors into $\phi$.  However, the two algorithms often work
together in concert -- each algorithm reducing the error introduced by 
the other.  

It is well-known that the reinitialization procedure can cause the zero 
level set of $\phi$ to shift from its ``true'' 
position~\cite{osher_book,seth99}.  Fortuitously, the process of 
enforcing the volume fraction constraint helps to reduce the magnitude of 
this shift by disallowing large shifts in the zero level set.  On the 
other hand, the volume fraction constraint algorithm can cause
a significant shift of the zero level set (especially when generating
initial conditions of a specified volume fraction).  If a narrow banding type 
of procedure is used, the zero level set could end up in a region where 
$\phi$ is not very smooth.  Reinitializing in a sufficiently wide band
prior to enforcing the volume fraction constraint helps to avoid this problem.
Reinitialization after enforcing the volume fraction constraint also helps
to keep the zero level set ``centered'' within the narrow band that is 
actively updated by the optimization algorithm.

\subsubsection{Stopping criteria for optimization loop}
Because our goal is to minimize $\cA$, the change in the interfacial 
surface area between iterations, $\Delta \cA_n$ is used as the stopping 
criteria.  Thus, the embedding function, $\phi$, is updated until the 
$\Delta \cA_n$ falls below a prescribe tolerance.

\subsubsection{Evaluation of surface integrals}
We numerically evaluate surface integrals by first converting them to volume
integrals via (\ref{eq:integral_conversion}) and then replacing the delta 
functions with smoothed approximates.  For example, using 
(\ref{eq:integral_conversion}), the Lagrange multiplier 
(\ref{eq:Lagrange_multiplier}) may be written in the form
\begin{equation}
\label{eq:curvature}
\lambda = -\frac{\int_\Omega (\nabla \cdot \bn) \delta(\phi) \left\| \nabla 
\phi \right\| dV}{\int_\Omega \delta(\phi) \left\| \nabla \phi \right\| dV}.
\end{equation}
In our computations, we use the following approximation for the
$\delta$-function~\cite{osher_book}:
\beq
\label{eq:delta}
\delta_\epsilon(\phi)\ =\
  \begin{cases}
    0 & \qquad \|\phi\|\ >\ \epsilon \\
    \frac{1}{2\epsilon}[1+\cos(\frac{\pi \|\phi\|}{\epsilon})] & 
    \qquad \|\phi\|\ \le\ \epsilon, 
  \end{cases}
\eeq
where $\epsilon$ is the width of the smoothed $\delta$-function.
Following the usual conventions for level set method calculations, we 
choose $\epsilon = 3 \max\{\Delta x_1,\Delta x_2, \Delta x_3\}$.

\subsubsection{Discretization of the mean curvature }
The discretization of the curvature term is performed by using the 
natural generalization of the first-order accurate scheme described by
Zhao, Chan, Merriman, and Osher~\cite{zhao1996}.

\subsubsection{Grid resolution requirements}
All of our numerical solutions were computed on a unit cube using uniform 
meshes with sizes between $100 \times 100 \times 100$ and 
$250 \times 250 \times 250$.  For relatively simple surfaces, such as the 
Schwartz P surface, the lowest resolution mesh was sufficient to obtain 
accurate values for the total surface area\footnote{The difference in the 
total surface area computed at the $100 \times 100 \times 100$ and
the $150 \times 150 \times 150$ resolutions was less than 
$10^{-4}$.}.  However, complex surfaces, where the distance between separate
sheets of the surface is smaller, required higher resolution meshes.  For
example, in order to obtain results comparable in accuracy to those for the 
Schwartz P and Schoen G surfaces with a mesh of size 
$200 \times 200 \times 200$, the Schwartz D surface calculations required
a grid of size $250 \times 250 \times 250$.
Finer meshes were also required to obtain accurate values for the mean
curvature and when finding optimal surfaces for volume fractions other
than $0.5$.

\subsubsection{Parallel computation}
Due to the high computational cost (both in time and memory) of the fully 
three-dimensional shape-optimization calculations, we implemented our 
algorithm using LSMLIB~\cite{LSMLIBWebSite2006}, a parallel level set method 
software library developed by one of the authors.  Rapid implementation of
our parallel simulation code and LSMLIB was achieved by leveraging the parallel
computing framework built into the Structured Adaptive Mesh Refinement 
Application Infrastructure 
(SAMRAI)~\cite{Hornung:ApplicationComplexity,SAMRAIWebSite2006} developed
at Lawrence Livermore National Laboratory (LLNL).

\section{Numerical Results}
\label{section5}
In this section, we present results obtained using the numerical optimization
algorithm developed in Section~\ref{section4}.  Simulations were run on  
Linux clusters\footnote{Runs were performed on either a 75 node dual Opteron 
cluster or a 128 node Intel P4 system.} using between 2 and 8 nodes depending 
on the size of the computation.  Results were visualized using the VisIt 
program developed at LLNL.

\subsection{Verification of numerical optimization algorithm}
\label{subsection5.1}
To verify our numerical optimization algorithm, we tested its ability to 
accurately compute a few structures whose local optimality are easily 
verified.  These results were computed on a grid of size 
$100 \times 100 \times 100$. 

\begin{figure}[htb!]
\centerline{
\includegraphics[width=1.75in]{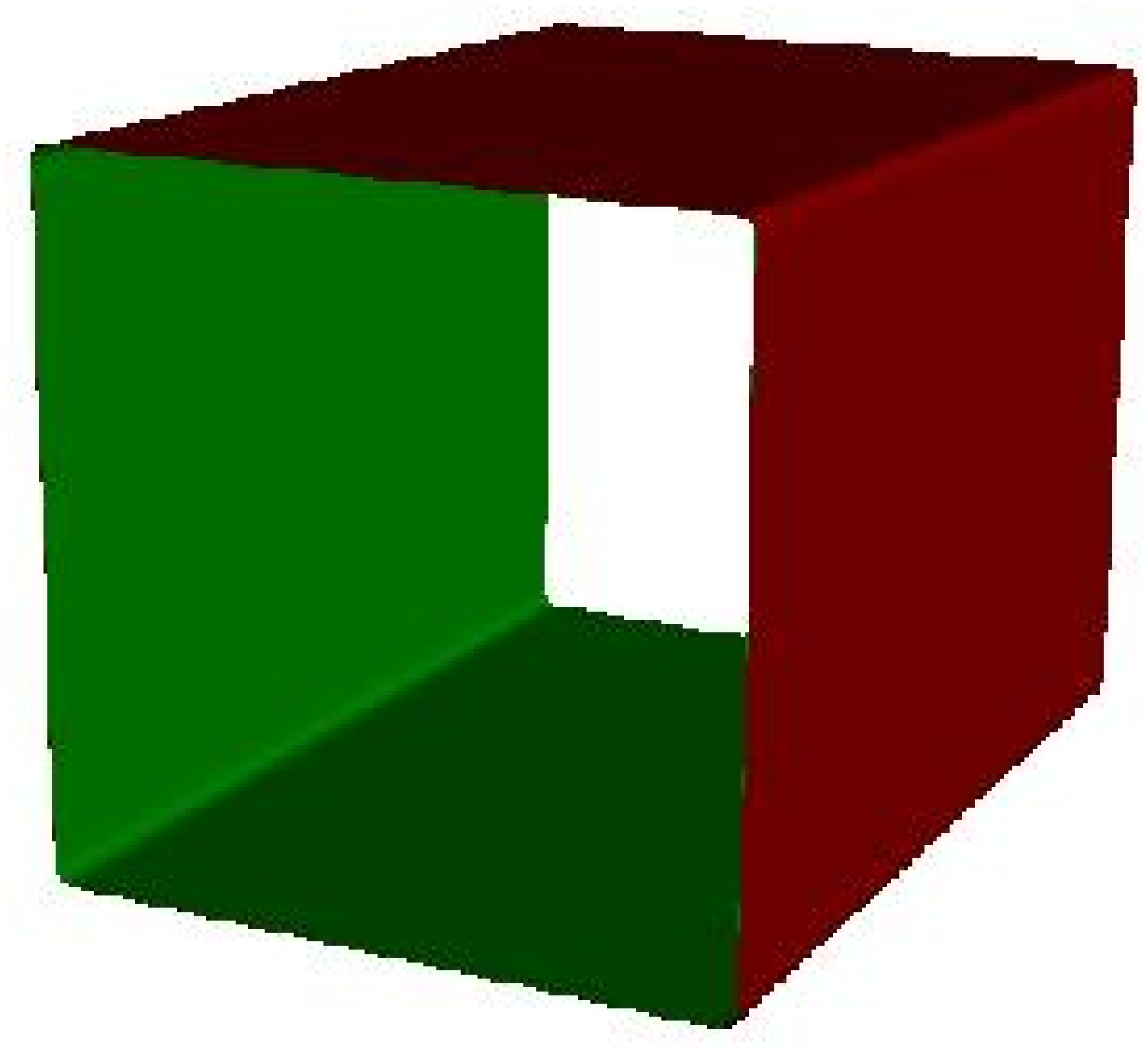} \hspace{0.5in}
\includegraphics[width=1.75in]{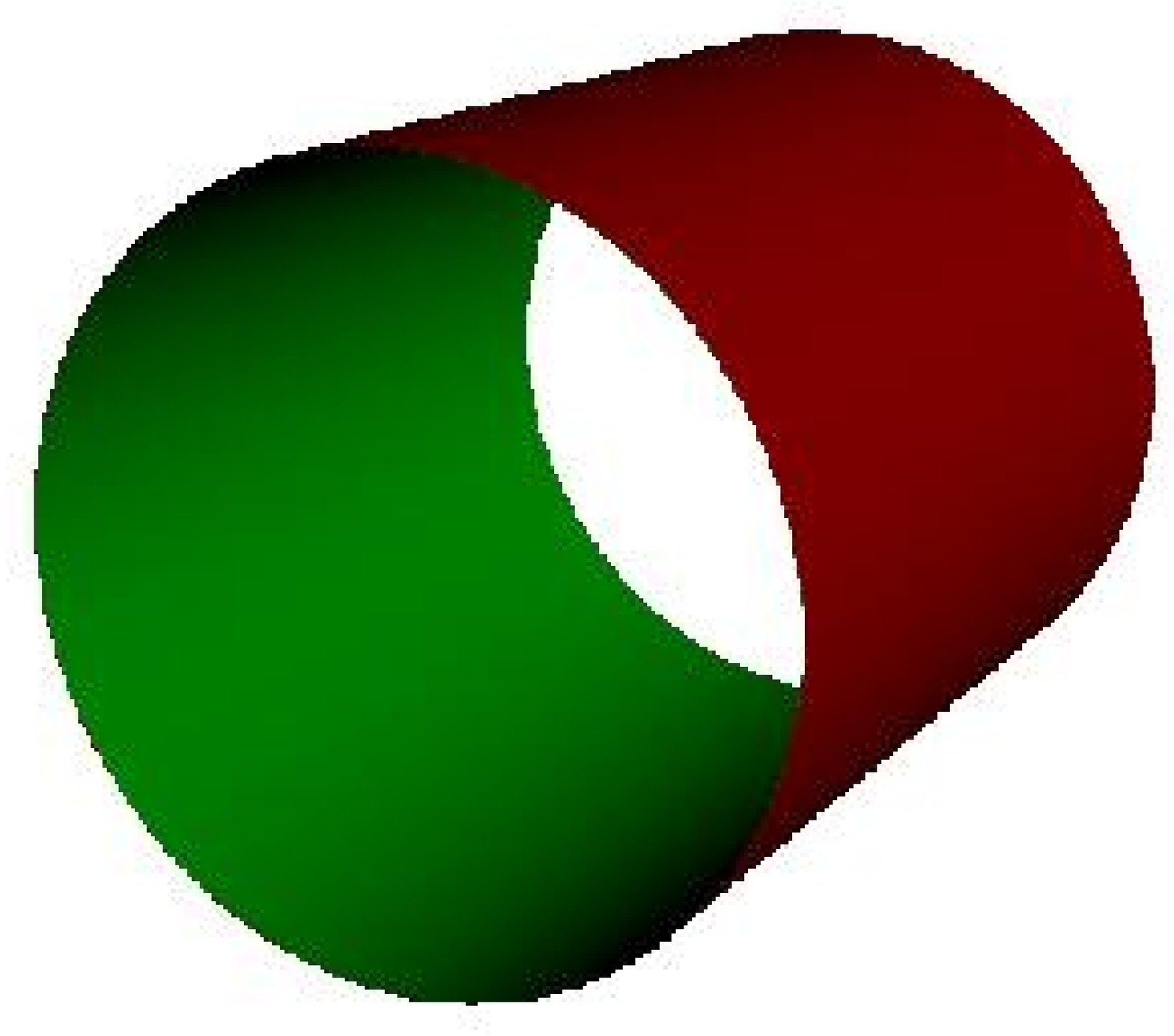} 
}
\caption{Evolution of an infinitely long square channel (left) to the
locally optimal cylindrical channel (right).}
\label{fig:algorithm_verification_cylinder}
\end{figure}

Our first test case is the cylinder.  Because it has constant mean curvature, 
Theorem~\ref{thm:extrema_total_surface_area} shows that it is a local
minimum of the total surface area under a volume fraction constraint.
From symmetry considerations, we expect that any infinitely long channel
should evolve towards the locally optimal cylinder structure. 
In our numerical calculations, we started the optimization algorithm with 
an infinitely long square channel as the initial surface. 
As expected, the square channel evolves to the optimal cylindrical channel
(see Figure~\ref{fig:algorithm_verification_cylinder}).

\begin{figure}[htb!]
\centerline{
\includegraphics[width=1.75in]{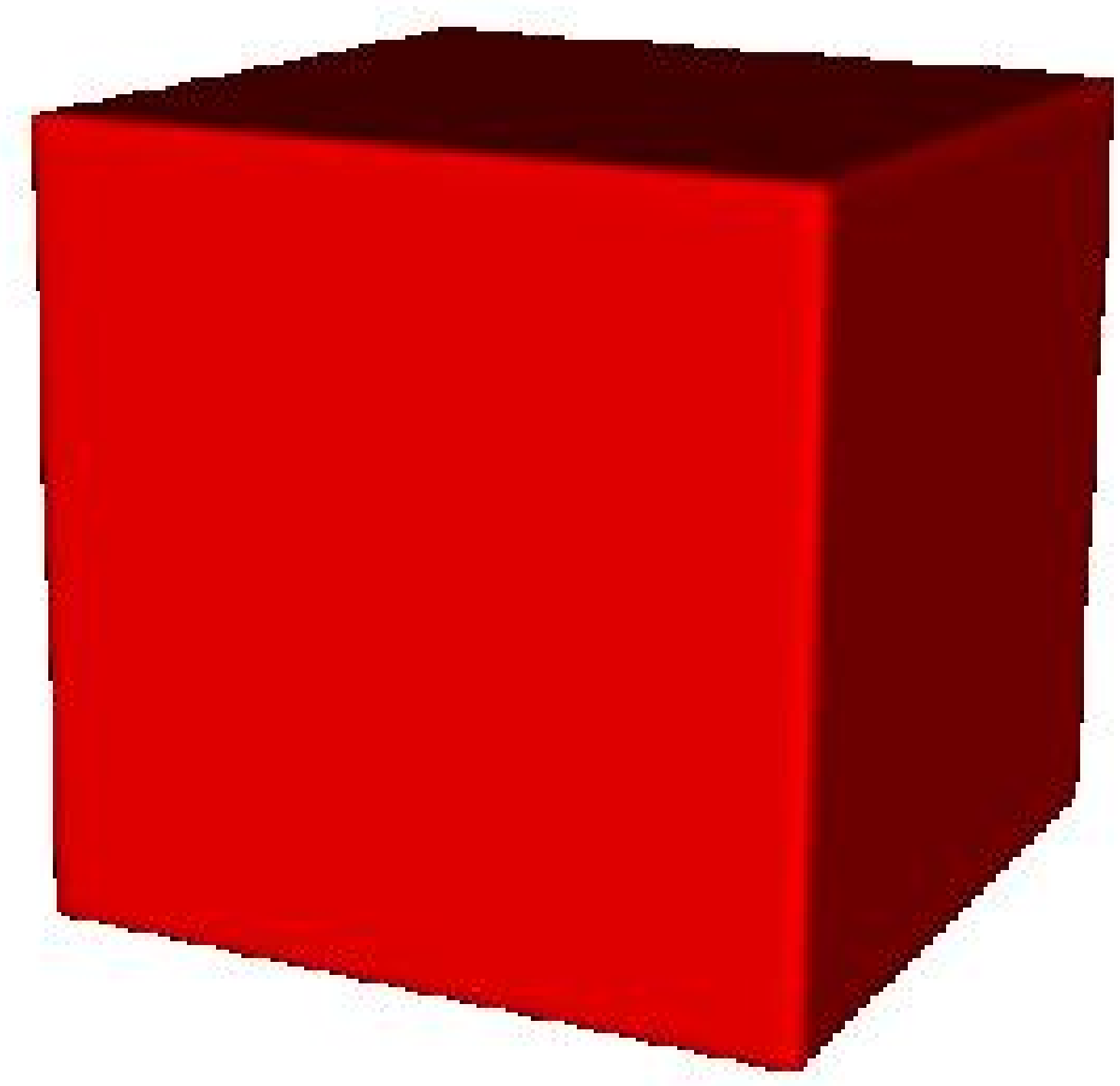} \hspace{0.5in}
\includegraphics[width=1.75in]{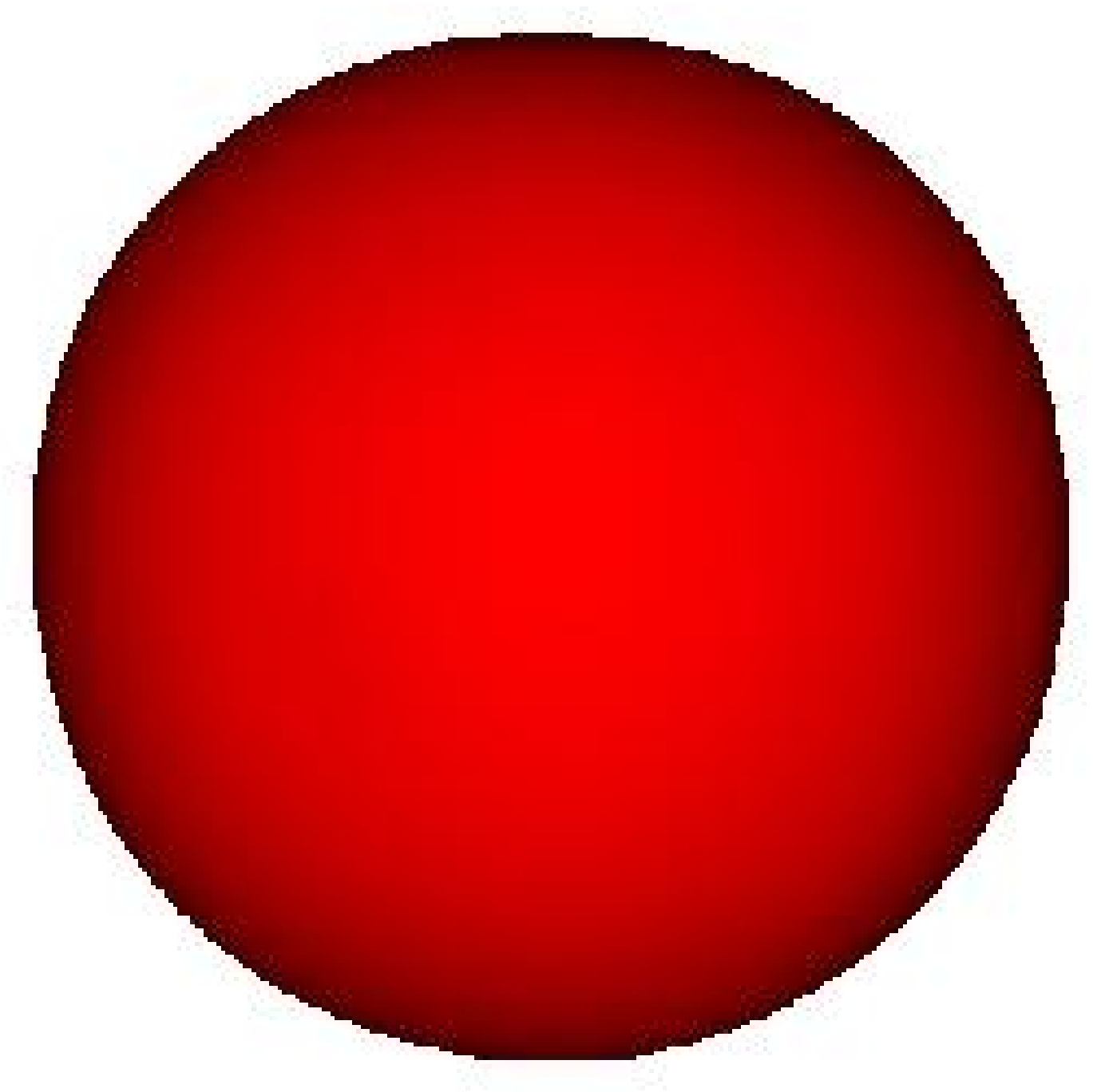}
}
\caption{Evolution of a cube surface (left) to the locally optimal sphere 
surface (right).}
\label{fig:algorithm_verification_sphere}
\end{figure}

Our second test case is the sphere.  Again, 
Theorem~\ref{thm:extrema_total_surface_area} shows it is a locally optimal
surface.  In this test, we started the optimization algorithm with 
an closed cube (chosen from symmetry considerations) as the initial surface. 
Figure~\ref{fig:algorithm_verification_sphere} shows that the cube evolves to the 
sphere, as expected. 

In both calculations, the volume fraction constraint algorithm was used to 
enforce $\volfrac_o = 0.5$.  However, the volume fraction did not fluctuate
much during the numerical simulations, so the volume fraction constraint 
algorithm was rarely invoked.  Also, notice that for both of these test cases, 
the level set calculations do not appear to have any difficulty handling sharp 
edges/corners.

\subsection{Local minimality of Schwartz P, Schwartz D, and 
            Schoen G surface areas}
\label{subsection5.2}
There are many known ways to compute the Schwartz P, Schwartz D, and Schoen G 
triply periodic minimal surfaces.  They can be characterized exactly using an 
Enneper-Weierstrass (complex integration) representation~\cite{Nits89},
generated as the local minima of the scalar order parameter Landau-Ginzburg
functional used to describe ordering phenomena in microemulsions~\cite{Gozd96a},
and approximated by Fourier series using the periodic nodal surface (PNS) 
expansion~\cite{Schn91,Mack93,Schw99,Gand01}.  It is worth pointing out
that the approach of minimizing the Landau-Ginzburg functional is actually 
a phase-field version of our approach for an appropriately chosen form of
the energy functional.  

Because minimal surfaces have zero mean curvature everywhere, 
Theorem~\ref{thm:extrema_total_surface_area} guarantees that the Schwartz P, 
Schwartz D, and Schoen G surfaces are locally optimal surfaces.  However,
it does not indicate whether these surfaces are local maxima or minima or 
saddle points of the total surface area.  The numerical optimization 
procedure developed in Section~\ref{subsection4.3} provides a means to answer 
this question.
Since the optimization procedure seeks surfaces with minimal total surface
area, locally minimal surfaces should be stable under 
perturbations.  In contrast, local maxima or saddle points of the total 
surface area should be unstable\footnote{Note that instability is a necessary
by not sufficient condition for a local maximum or saddle point while
stability is a sufficient by not necessary condition for a local minimum.  
Thus, we can only safely draw conclusions when the surface is a local minimum
of the total surface area.}.  
\begin{figure}[htb!]
\centerline{
\includegraphics[width=2.0in]{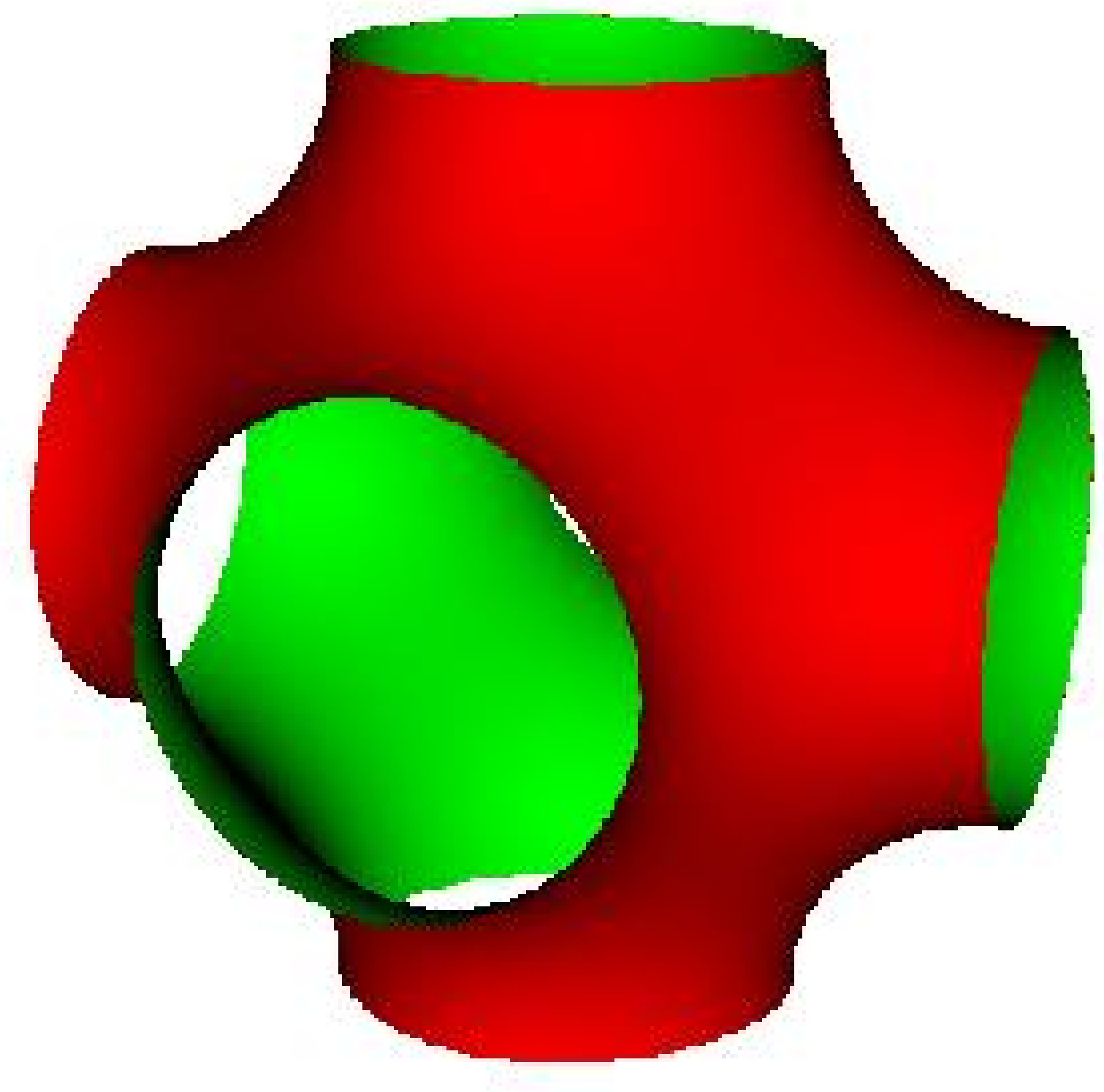} \hspace{0.0in}
\includegraphics[width=2.0in]{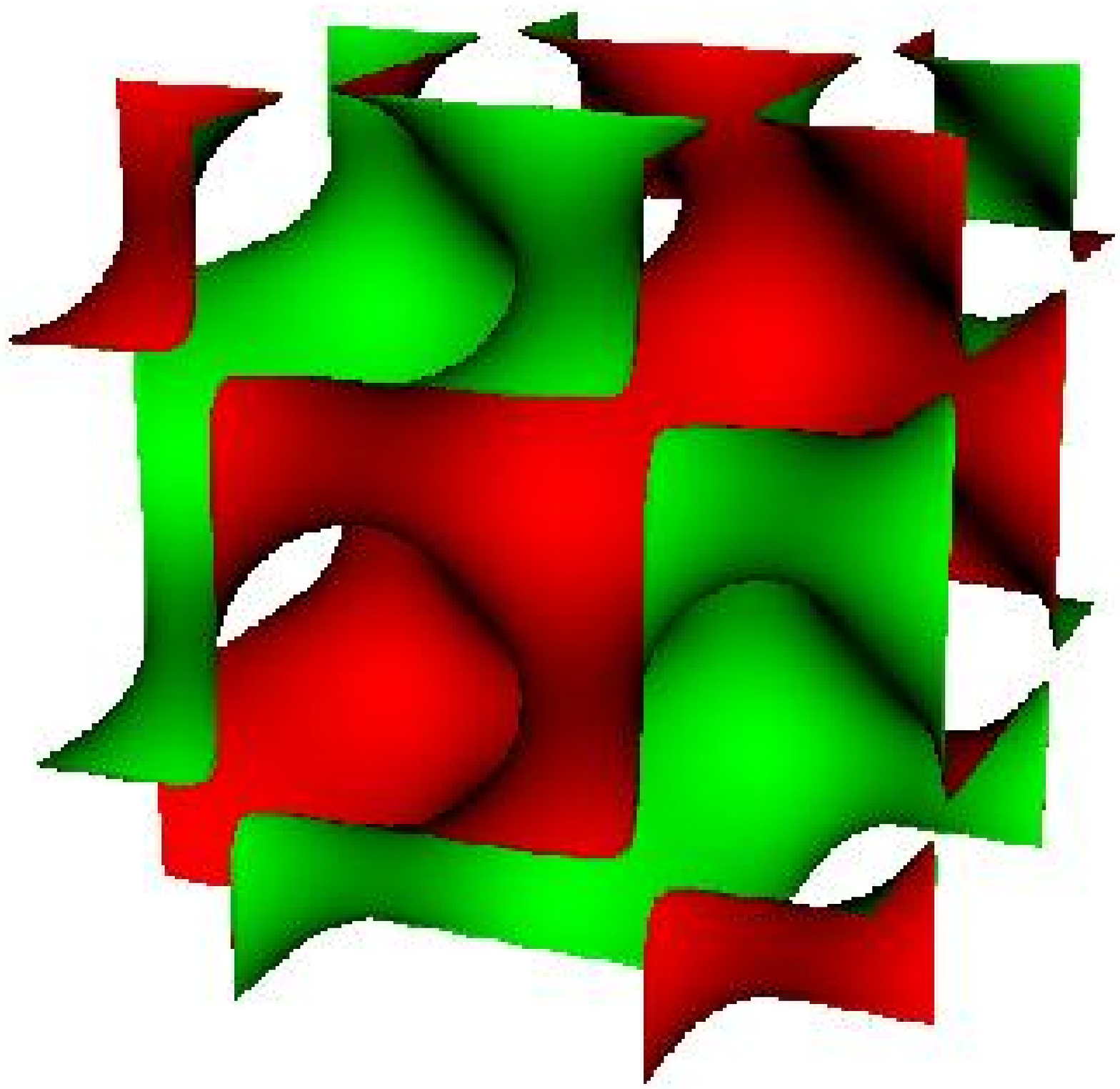} \hspace{0.0in}
\includegraphics[width=2.0in]{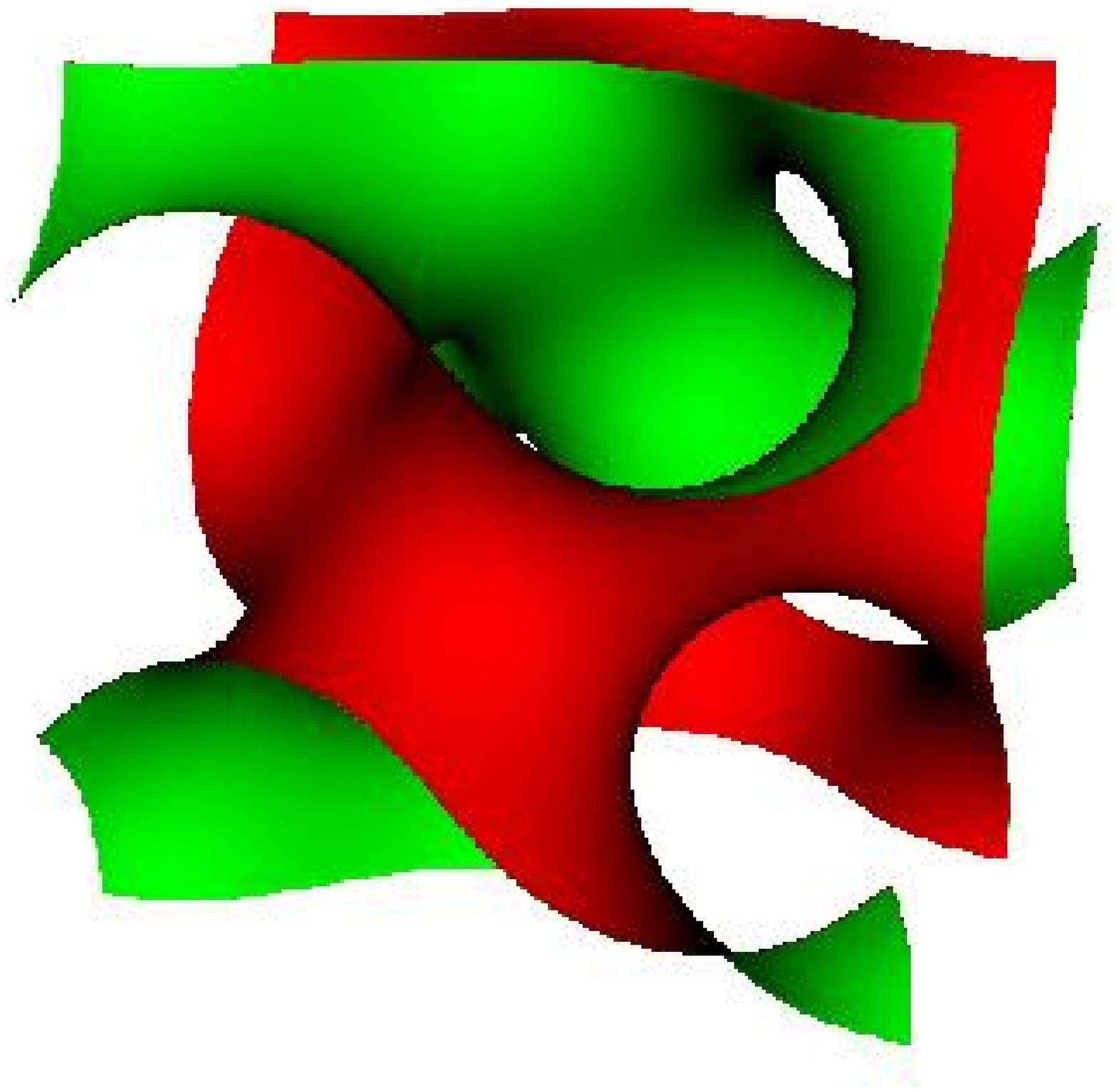}
}
\caption{Locally optimal surfaces for $\volfrac_o=1/2$ starting 
from the
PNS approximation to the Schwartz P surface (left), Schwartz D surface 
(middle), Schoen G surface (right).  Notice that the final surfaces are
precisely the Schwartz P, Schwartz D, and Schoen G surfaces, which indicates
that all three of these surfaces are local minima of the total surface
area when the volume fraction is $1/2$.  These results were computed using 
a grid of size $150 \times 150 \times 150$.}
\label{fig:minimality_of_common_surfaces}
\end{figure}

In our numerical simulations, we start the optimization procedure using
the PNS approximations\footnote{In the PNS approximation, each term consists 
of sine and cosine functions where the high-order terms have a high reciprocal 
lattice vector norm \cite{Schw99}.} (which are very good) as the initial 
$\phi$.  As can be seen in Figure~\ref{fig:minimality_of_common_surfaces}, 
the optimization procedure evolves to the correct triply periodic minimal
surface.  Therefore, we can conclude that the Schwartz P, Schwartz D, and 
Schoen G surfaces are all local minima of total surface area when 
$\volfrac_o=1/2$.

\subsubsection{Width of basins of attraction for Schwartz P, Schwartz D
and Schoen G}
\label{subsubsection5.2.1}
By starting the optimization procedure with other initial configurations
of carefully controlled symmetry, we were able to qualitatively explore the 
width of the basin of attraction for the Schwartz P, Schwartz D, and Schoen 
G surfaces.

\begin{figure}[htb!]
\centerline{
\includegraphics[width=1.8in]{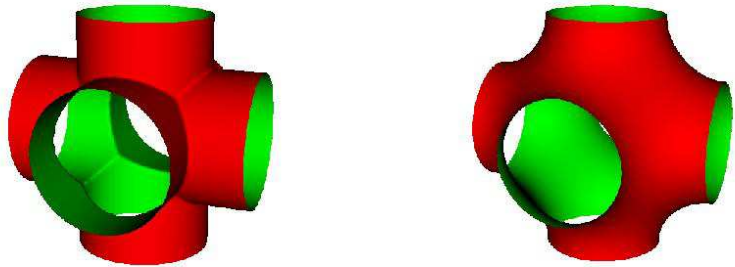} \hspace{0.5in}
\includegraphics[width=1.8in]{figs/P_0_5_final}
}
\caption{Evolution of triply periodic circular channels (left) to the
Schwartz P surface (right).}
\label{fig:triply_circ_chan}
\end{figure}

\begin{figure}[htb!]
\centerline{
\includegraphics[width=1.8in]{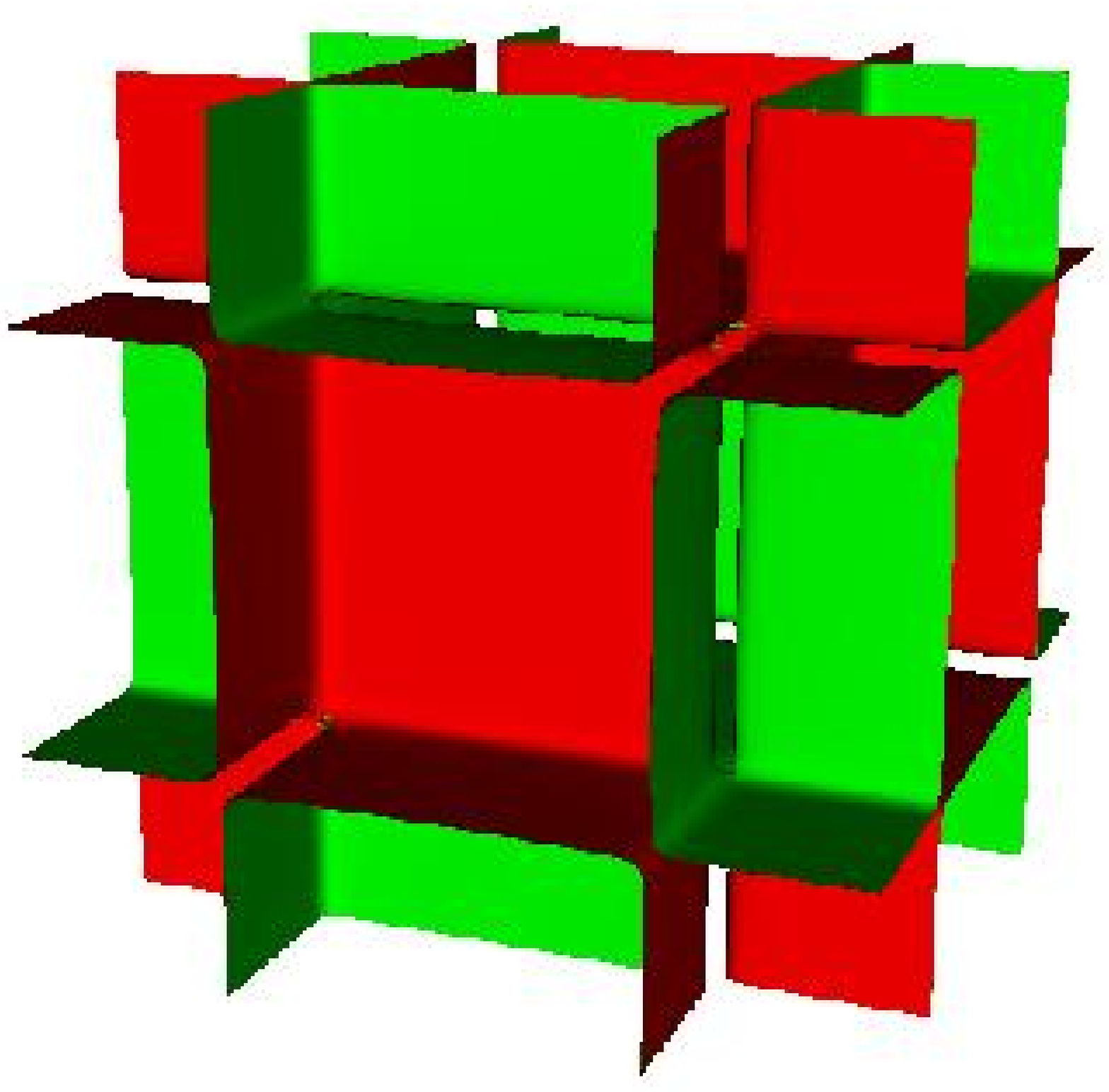} 
\hspace{0.5in}
\includegraphics[width=1.8in]{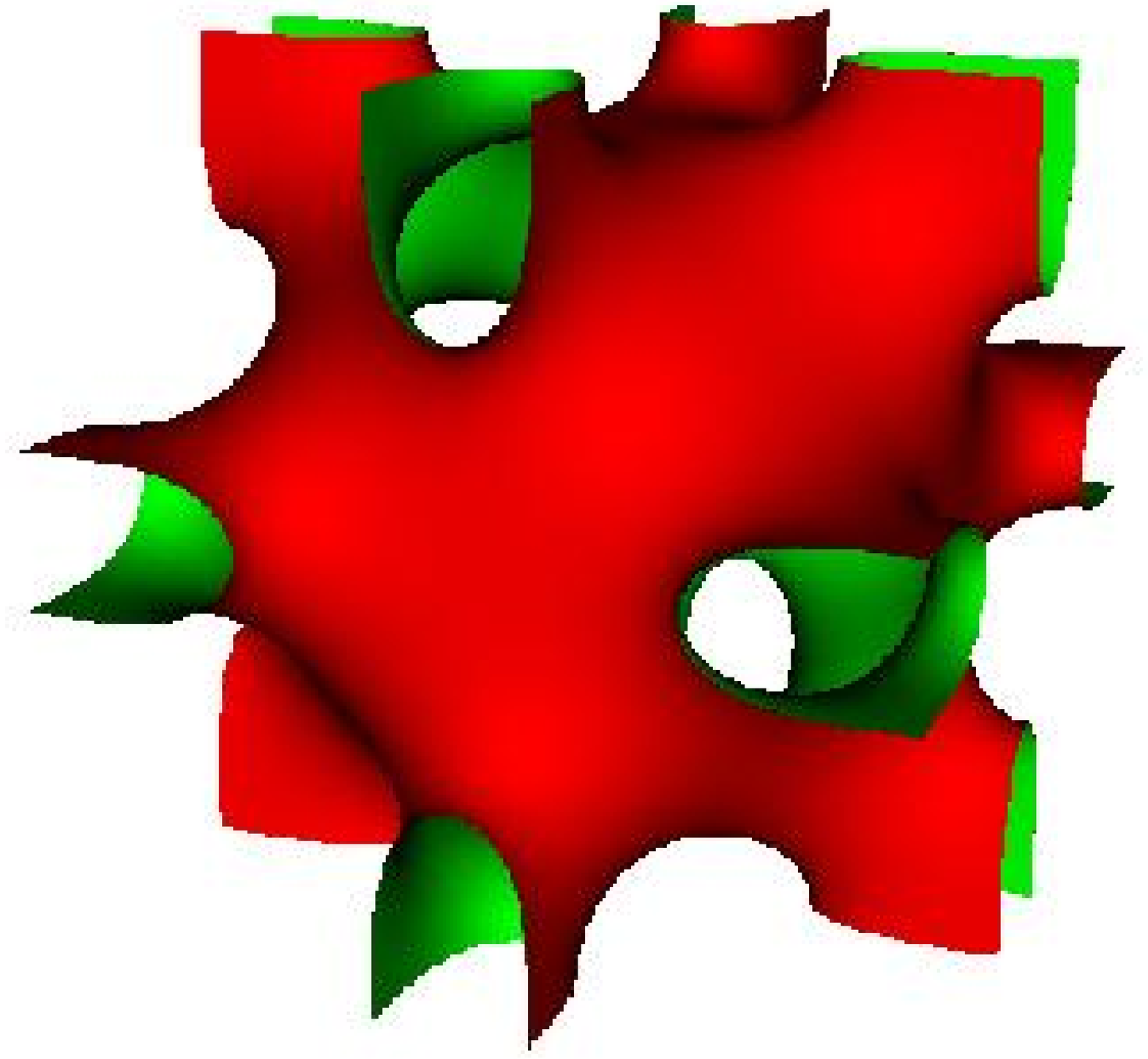}
}
\caption{Evolution of a triply periodic surface with the symmetry of the 
Schwartz P surface to a new high genus minimal surface (right).  The initial 
structure (left) is obtained by modifying the PNS approximation to the 
Schwartz P surface so that the second term dominates.} 
\label{fig:GBMK_P_second}
\end{figure}

For the Schwartz P surface, we considered two different initial configurations:
(1) triply periodic circular channels (see Figure~\ref{fig:triply_circ_chan})
and 
(2) $\phi$ set to a function dominated by the second term in the PNS 
approximation of the Schwartz P surface. 
As seen in Figure~\ref{fig:triply_circ_chan}, the triply periodic circular
channels configurations evolves to the Schwartz P surface.  However, the 
second initial configuration leads to a new high genus\footnote{The genus 
describes how many holes are in a closed surface and therefore is an 
integer number.  In the infinitely periodic case, we describe the genus in 
a unit cell, which will be a finite number~\cite{An90}.} surface with the 
same symmetry as the Schwartz P surface (see Figure~\ref{fig:GBMK_P_second}). 
The symmetry of the computed locally optimal solution is expected from 
symmetry of the initial configuration, but the genus of the solution is 
a bit unexpected. Clearly, this result is related to the width of the basin 
of attraction for the Schwartz P surface.

\begin{figure}[htb!]
\centerline{
\includegraphics[width=1.8in]{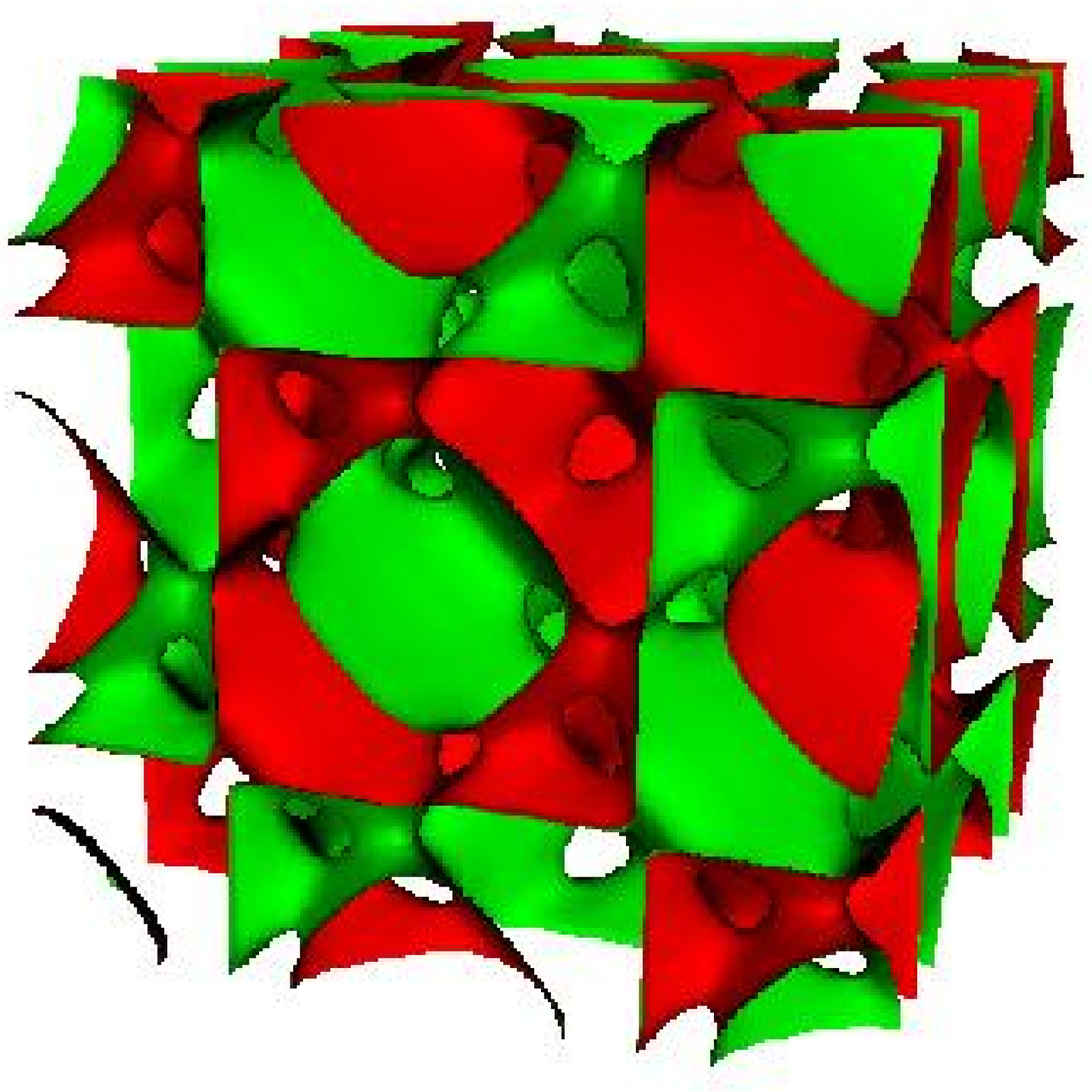} 
\hspace{0.5in}
\includegraphics[width=1.8in]{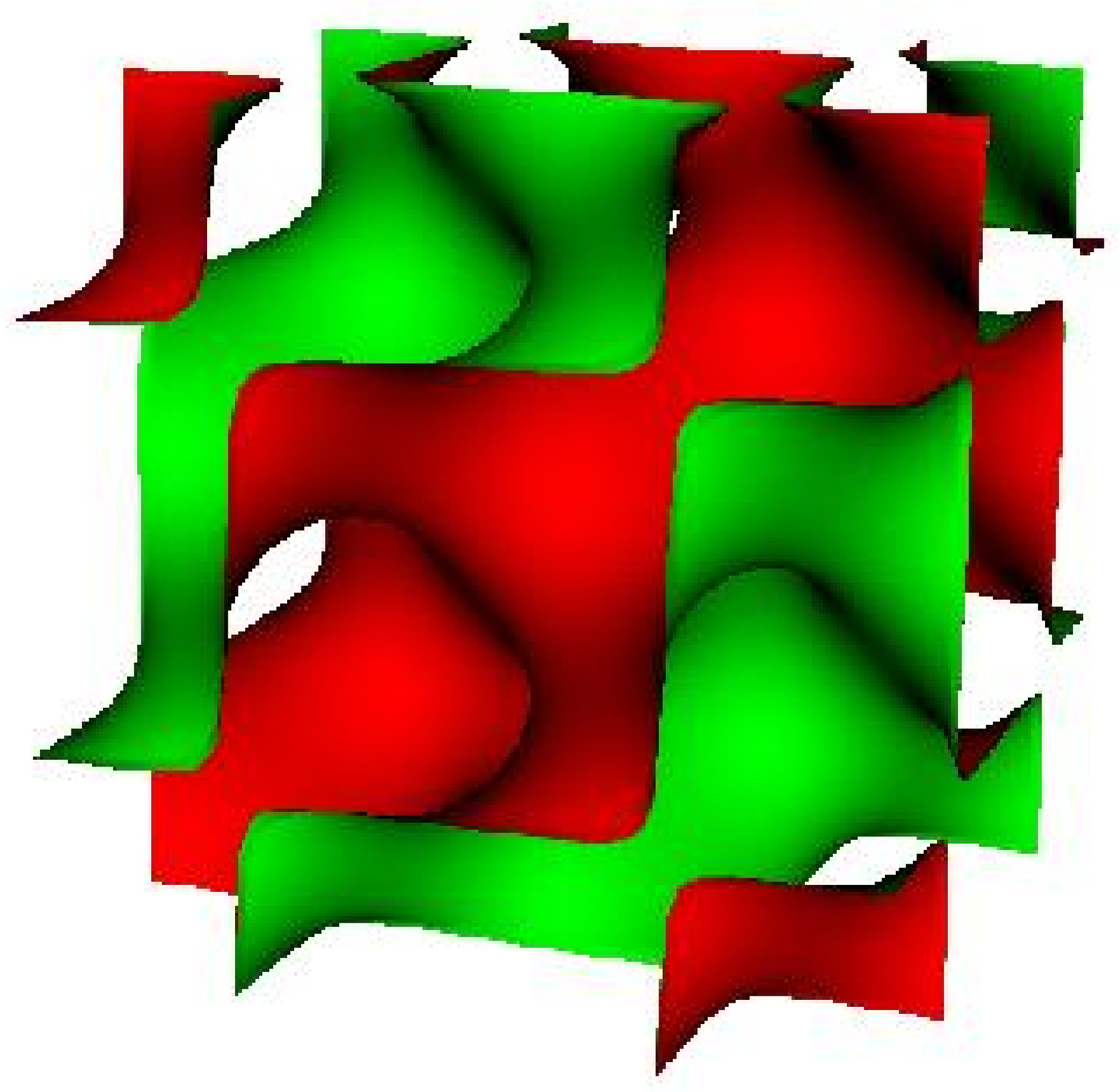}
}
\caption{Evolution of a triply periodic surface with the symmetry of the 
Schwartz D surface to a the Schwartz D surface (right).  The initial 
structure (left) is obtained by modifying the PNS approximation to the 
Schwartz D surface so that the second term dominates.} 
\label{fig:GBMK_D_second}
\end{figure}

\begin{figure}[htb!]
\centerline{
\includegraphics[width=1.8in]{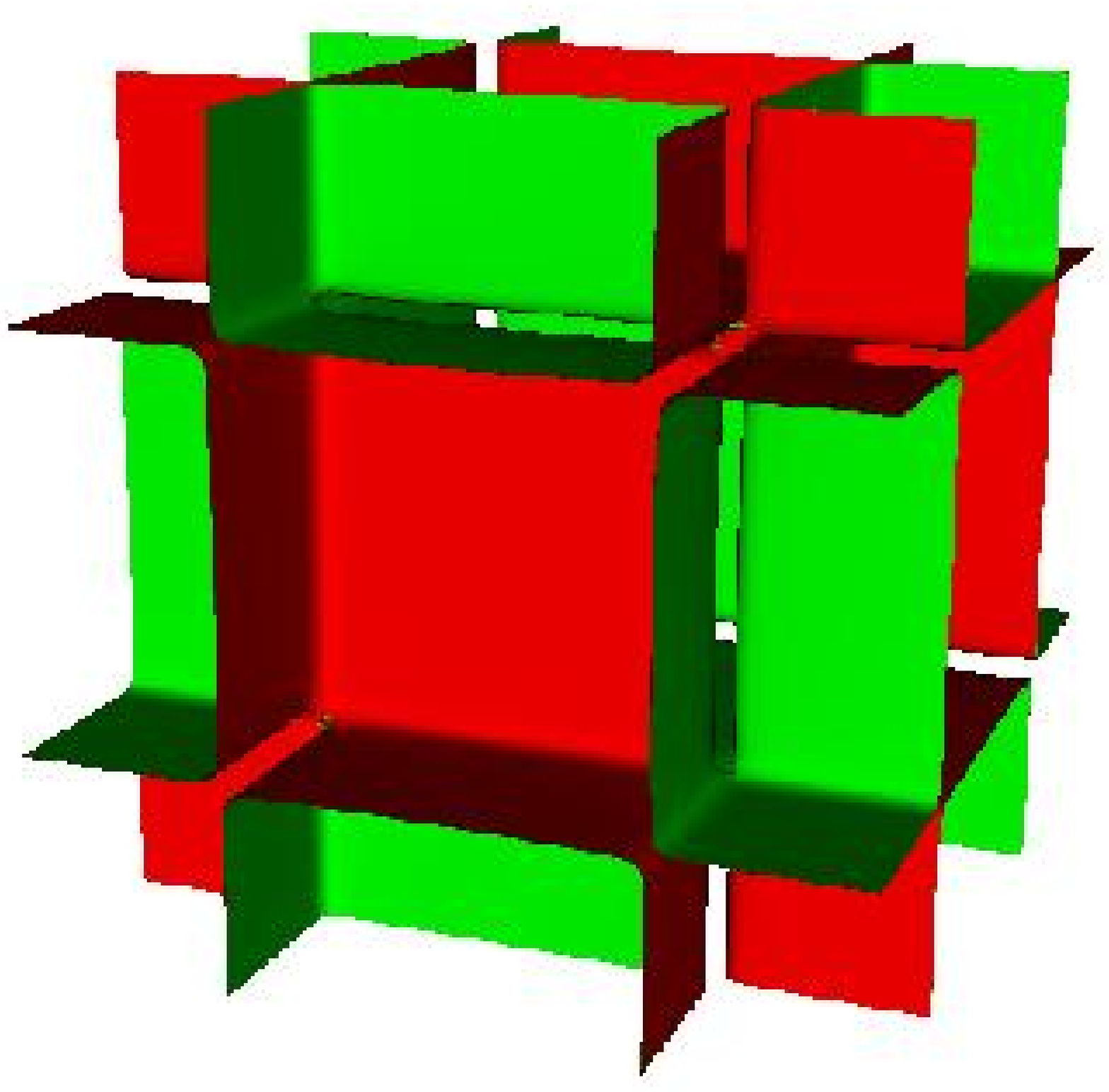} 
\hspace{0.5in}
\includegraphics[width=1.8in]{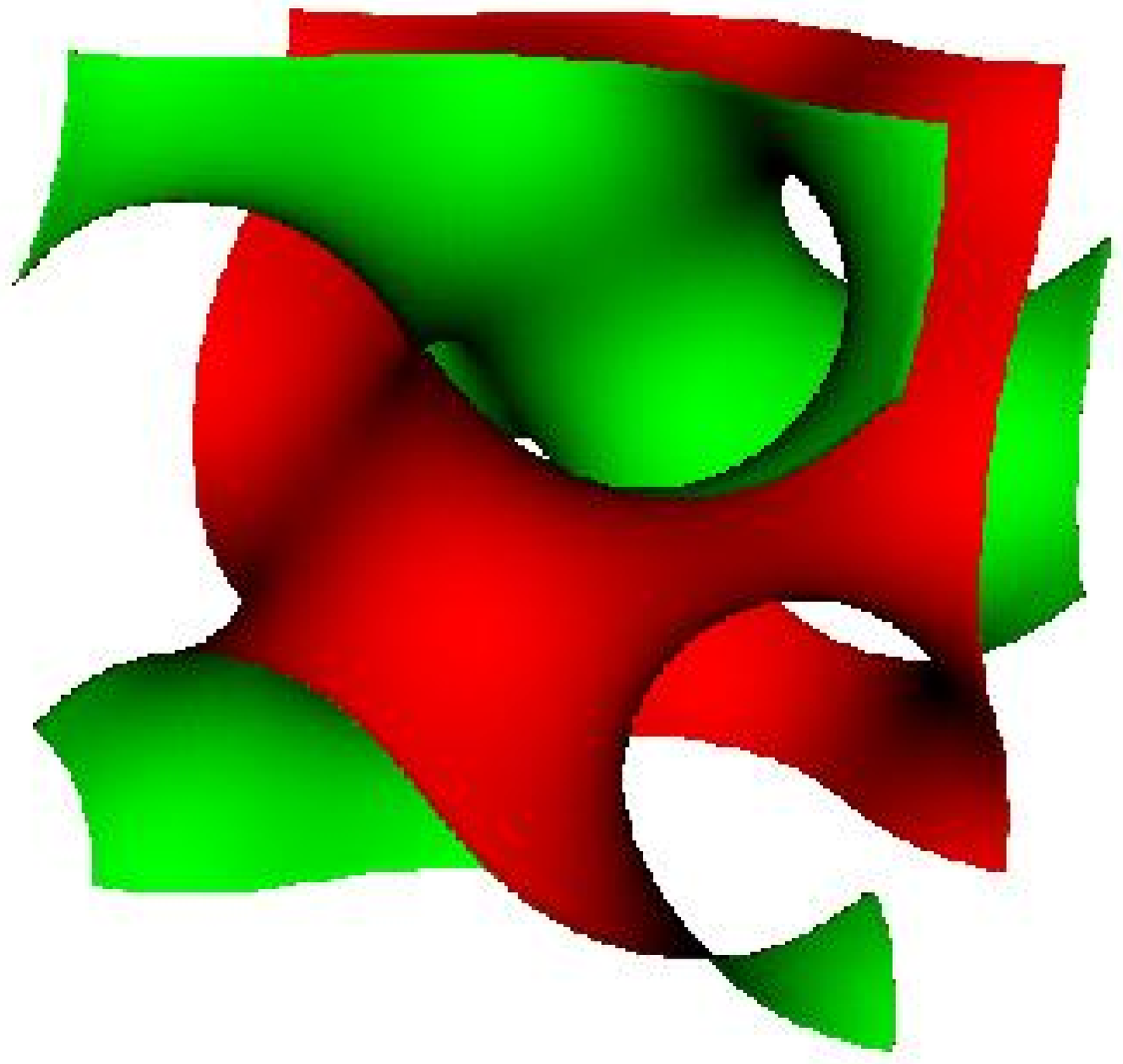}
}
\caption{Evolution of a triply periodic surface with the symmetry of the 
Schoen G surface to a the Schoen G surface (right).  The initial 
structure (left) is obtained by modifying the PNS approximation to the 
Schoen G surface so that the second term dominates.} 
\label{fig:GBMK_G_second}
\end{figure}

For the Schwartz D and Schoen G surfaces, we generated initial configurations
by setting the initial condition for $\phi$ to be a function dominated by 
the second term in the PNS approximation of the Schwartz D and Schoen G 
surfaces, respectively.
This procedure allowed us to preserve the symmetry of the of respective
surface, while significantly perturbing the actual surface.
Figure~\ref{fig:GBMK_D_second} shows that the Schwartz D surface is the
locally optimal structure even if the second term in the PNS approximation 
is dominant.  Even though the initial configuration in 
Figure~\ref{fig:GBMK_D_second} looks completely different from the Schwartz D 
surface, it converges to Schwartz D surface. This result implies that the 
Schwartz D has a relatively wide basin of attraction.  
The Schoen G surface shows a similar wide basin of attraction (see
Figure~\ref{fig:GBMK_G_second}). It should be noted that the three locally
optimal structures all have zero mean curvature at every point on their 
surfaces.

\subsection{Locally optimal structures when $\volfrac_o \neq 1/2$.}
\label{subsection5.3}
In this section, we examine three families of surfaces with the symmetry and 
topology as the Schwartz P, Schwartz D, and Schoen G surfaces
but with different values for the volume fraction of phase $1$.
We generated initial configurations with the desired symmetry, topology, and 
volume fractions by taking the PNS Schwartz P, Schwartz D and Schoen G 
approximations, respectively, and using the volume constraint algorithm to 
impose the desired volume constraint. We note that to obtain initial 
conditions for some of the volume fractions reported, we used continuation 
in the volume fraction to improve the convergence properties of the Newton 
iteration scheme.  These initial configurations were then allowed to evolve 
to the nearest optimal structure. 
The total surface area $\cA$ is computed using (\ref{eq:surface}) and 
the mean curvature $H$ is obtained from (\ref{eq:curvature}) using
the relation $H = 0.5\lambda$.  In both calculations, we have used the 
smoothed delta function, $\delta_\epsilon$, defined in (\ref{eq:delta}).

\subsubsection{Schwartz P surface family}
\label{subsubsection5.3.1}
\begin{figure}[htb!]
\centerline{
\includegraphics[width=1.6in]{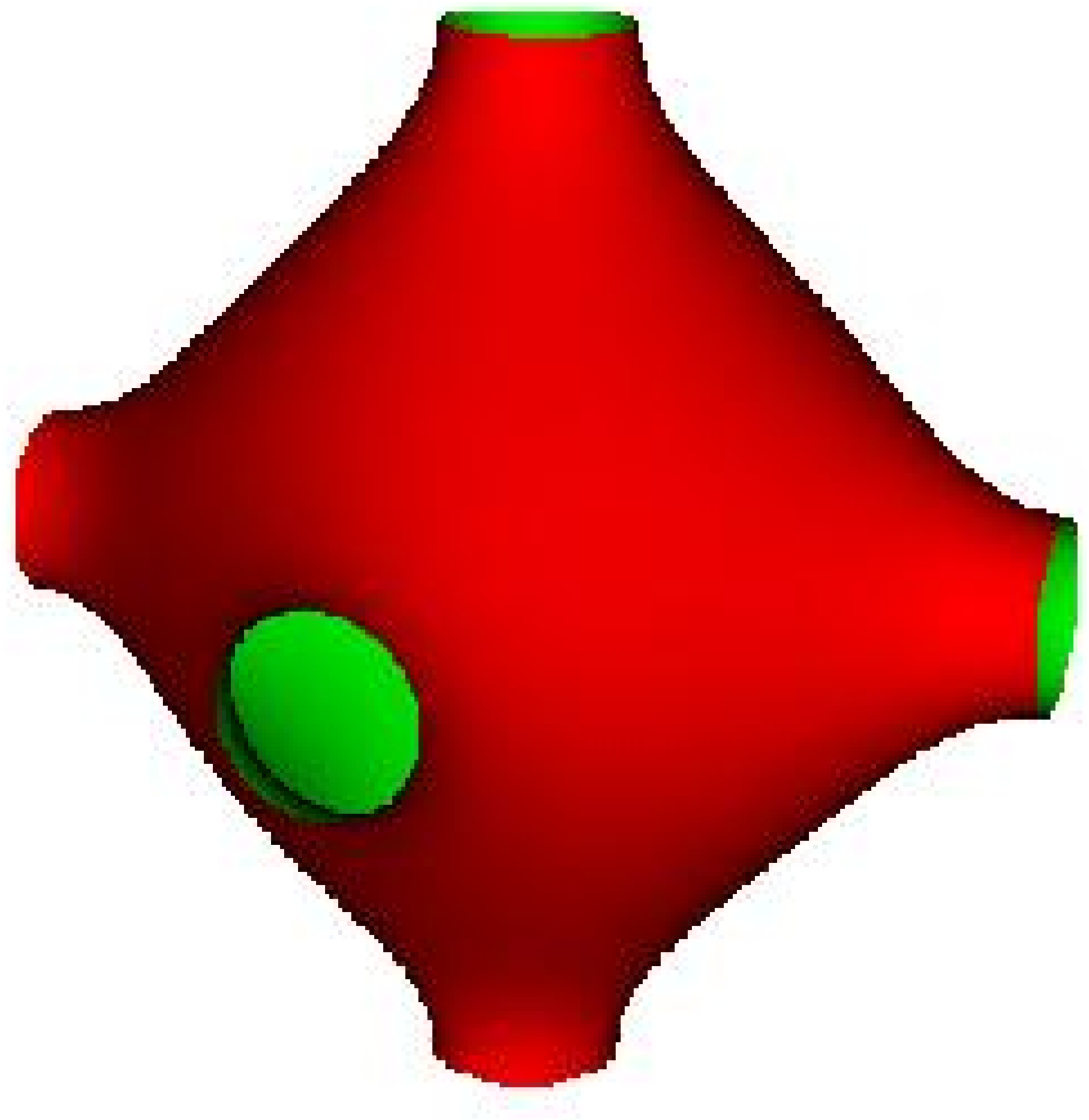} \hspace{0.0in}
\includegraphics[width=1.6in]{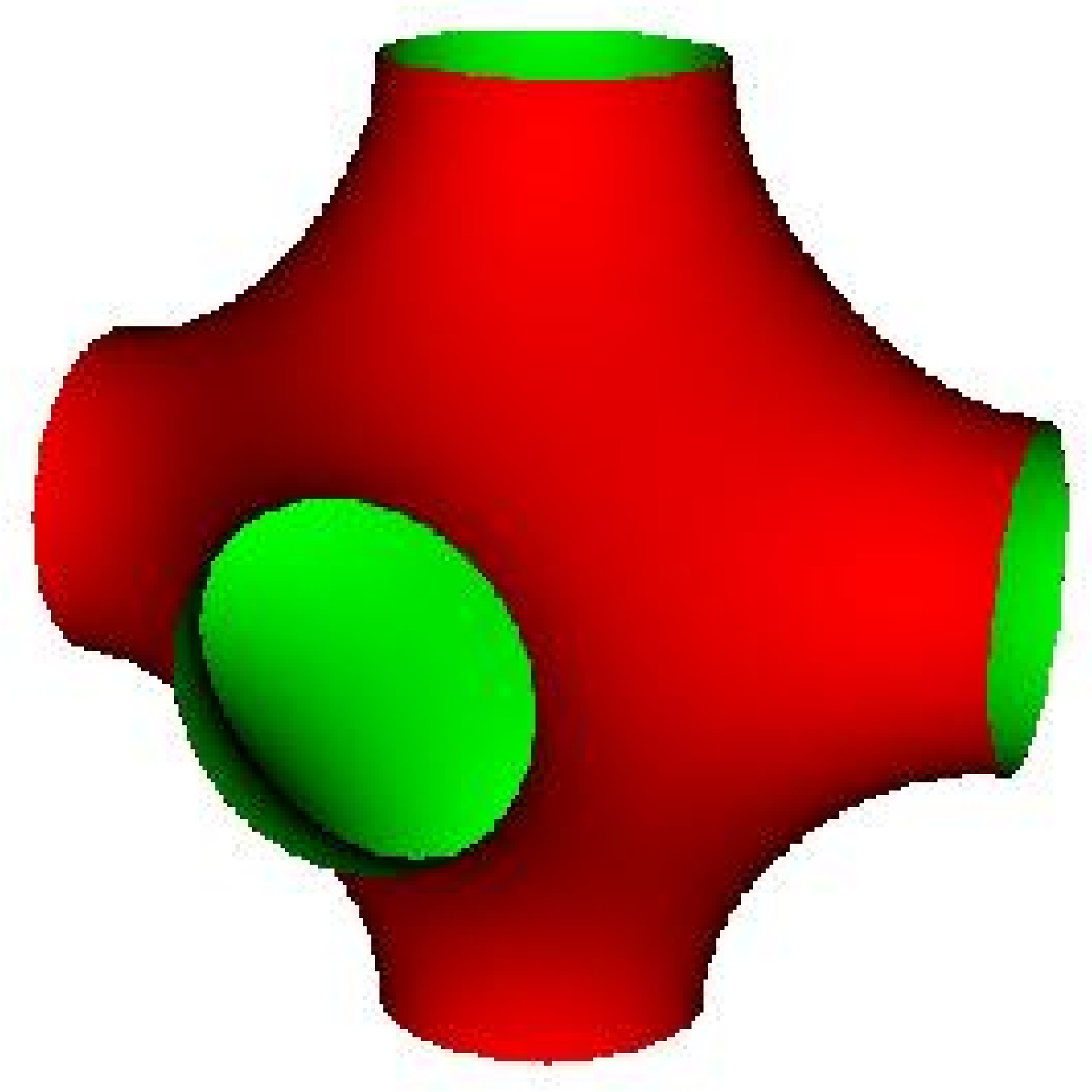} \hspace{0.0in}
\includegraphics[width=1.6in]{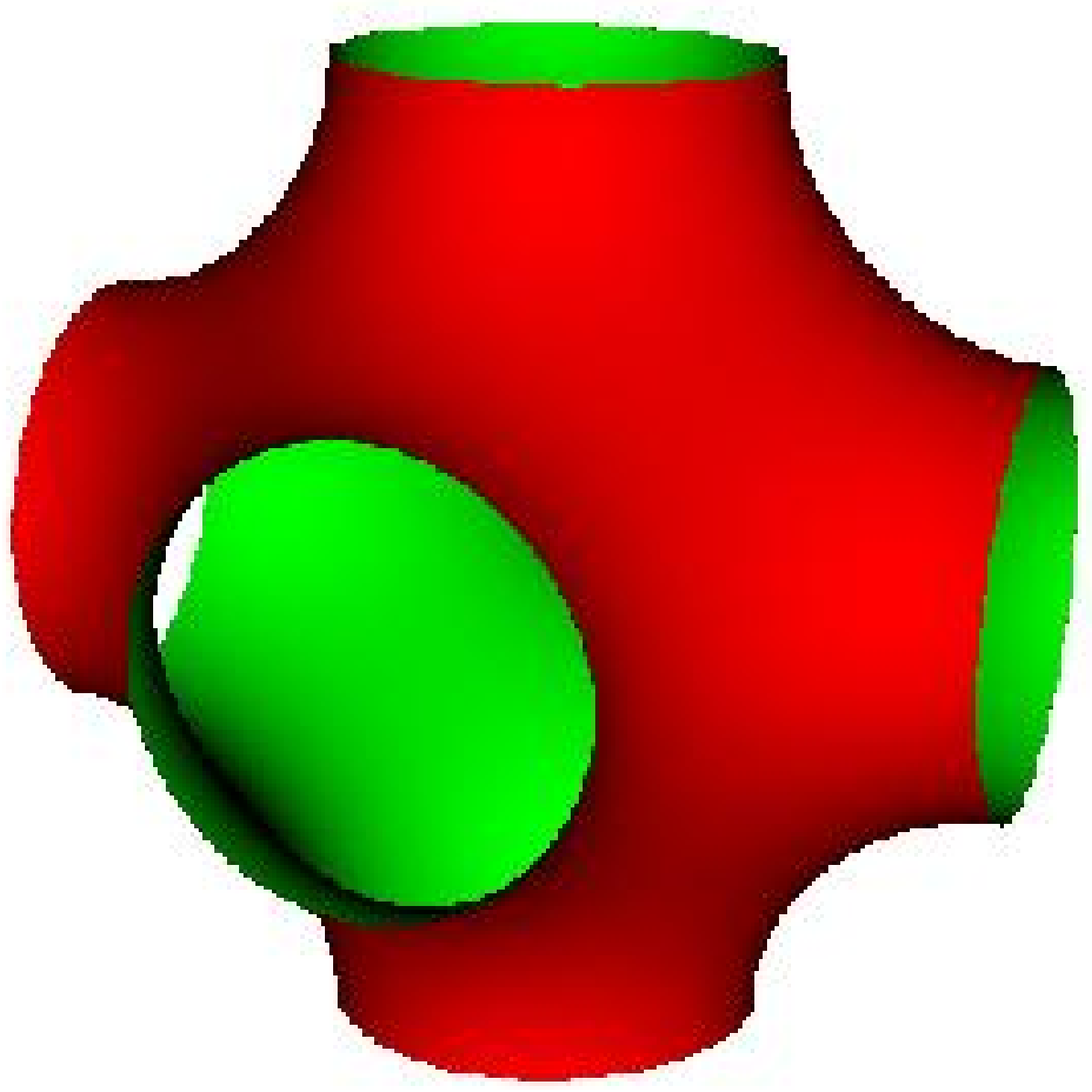}
}
\vspace{0.1in}
\centerline{
\includegraphics[width=1.6in]{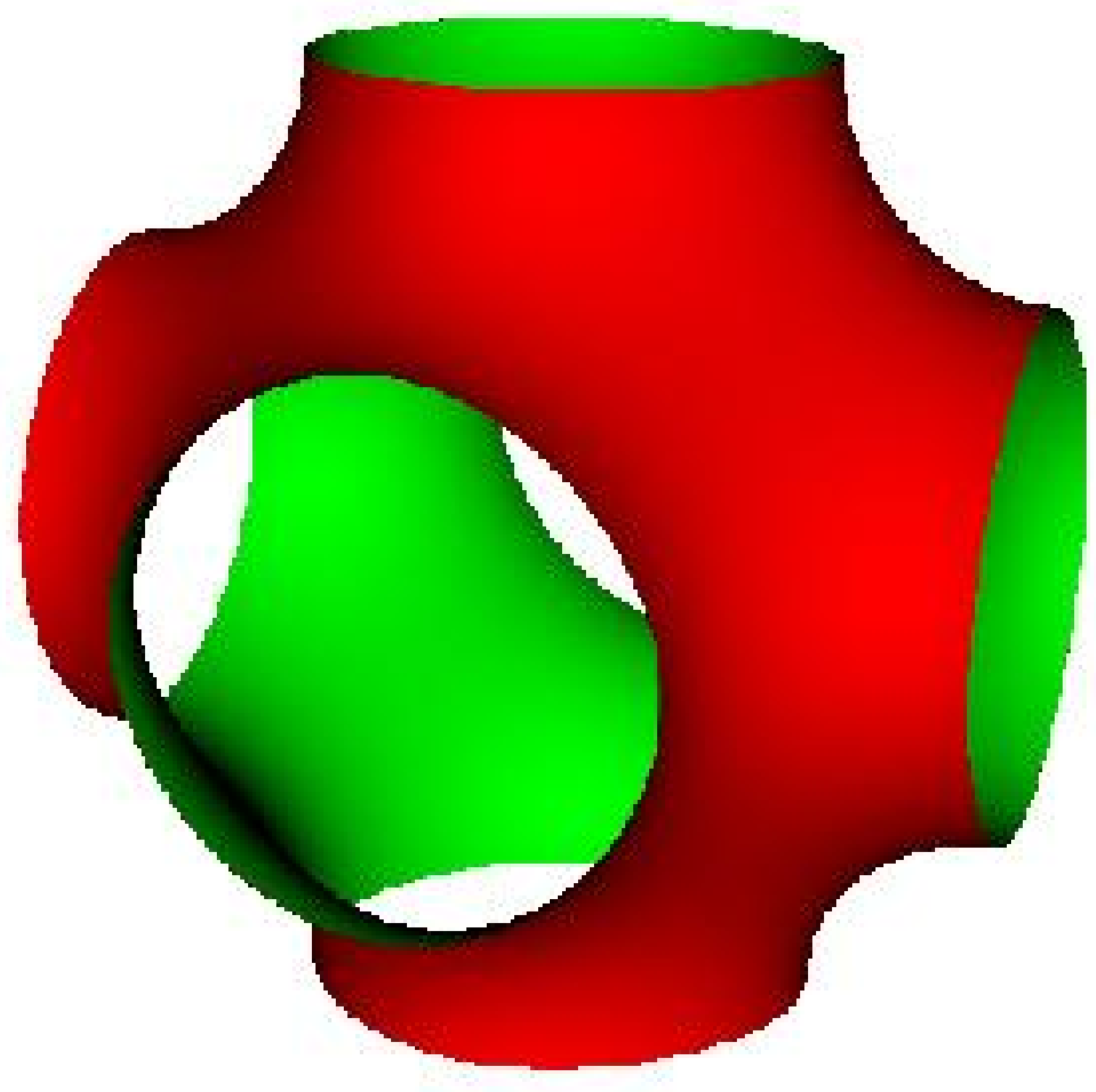} \hspace{0.0in}
\includegraphics[width=1.6in]{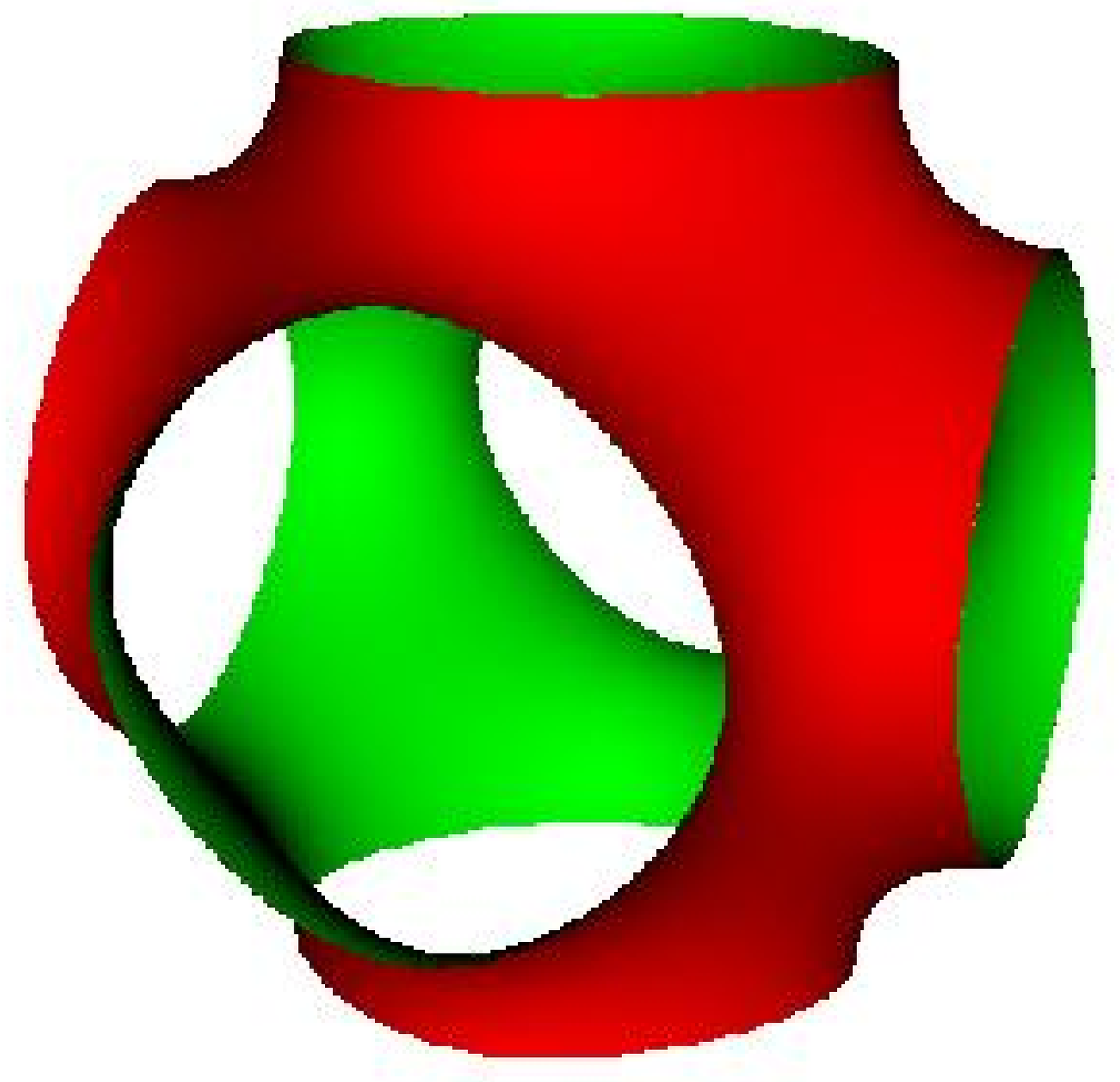} \hspace{0.0in}
\includegraphics[width=1.6in]{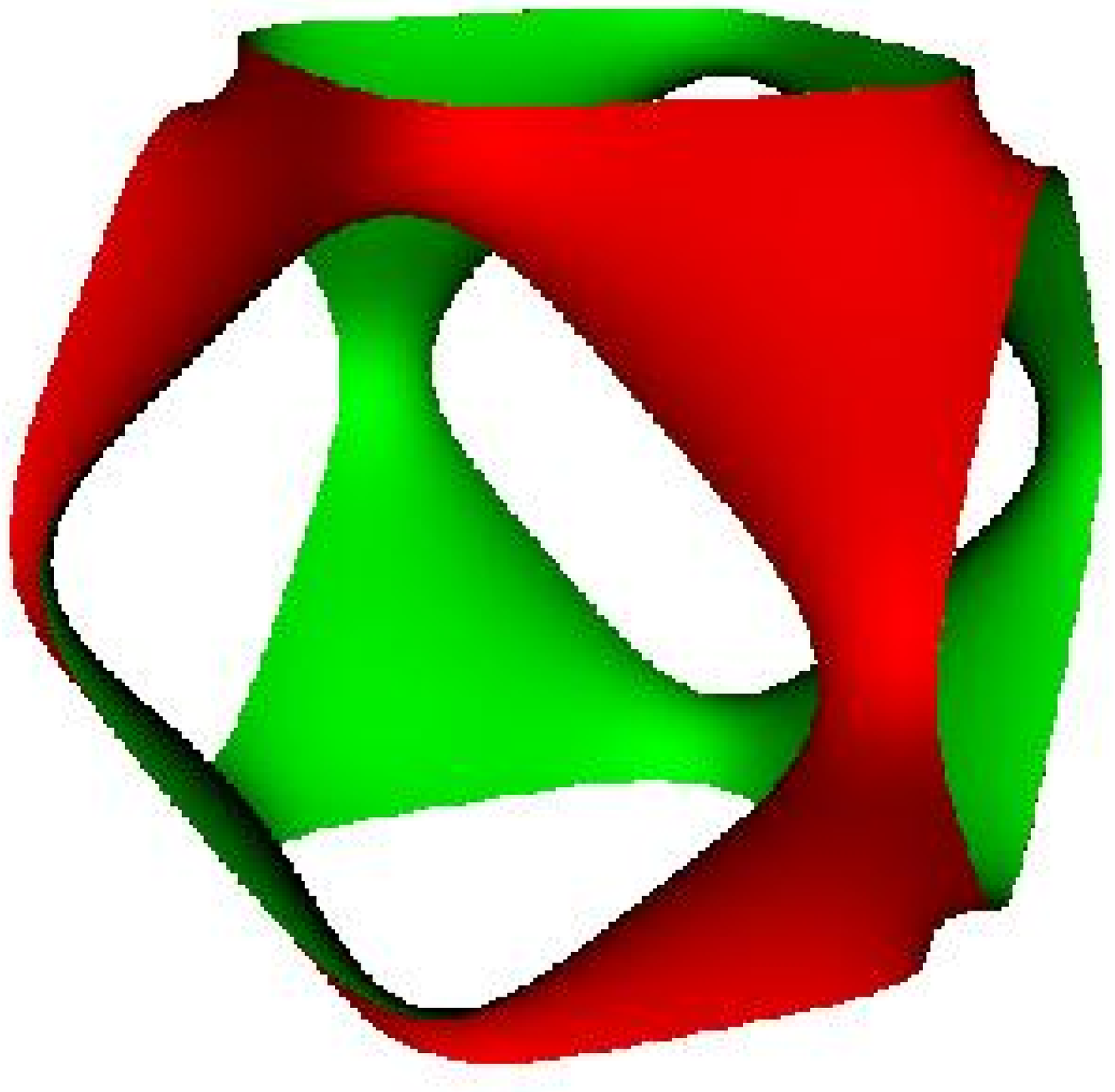}
}
\caption{Family of surfaces possessing the symmetry and topology
of the Schwartz P surface but with different volume fractions for phase 1. 
Upper left to upper right panel: $\volfrac_o = 0.25, 0.35, 0.45$. 
Lower left to lower right panel: $\volfrac_o = 0.55, 0.65, 0.75$.}
\label{fig:P_different_volume}
\end{figure}

We carried out the shape optimization procedure for the following volume 
fractions: 0.25, 0.3, 0.35, 0.4, 0.45, 0.5, 0.55, 0.6, 0.65, 0.7, 0.75.
Figure~\ref{fig:P_different_volume} shows the six different Schwartz P-type 
surfaces having specified volume fractions.
These calculations were performed on a grid of size $200 \times 200 
\times 200$. For a typical simulation, the area stopping criterion 
$\Delta \cA_n < 10^{-6}$ was reached after approximately 29 hours for a
4 node parallel calculations on the Opteron cluster.  The computation time 
in general increases for volume fractions further from $\volfrac_o = 0.5$ 
because the volume fraction constraint algorithm is invoked more often.

\begin{figure}[htb!]
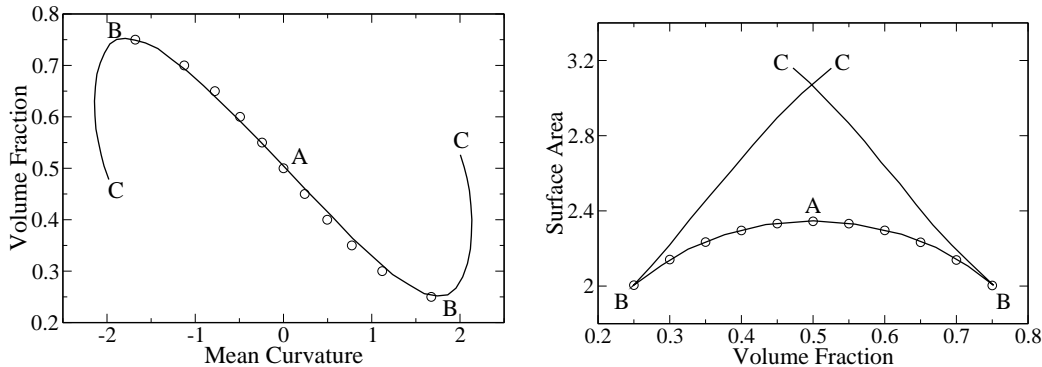

\centerline{
\includegraphics[width=2.6in]{figs/P_curv_vol_200} \hspace{0.1in}
\includegraphics[width=2.6in]{figs/P_vol_area_200}
}
\caption{Relationships between the volume fraction, mean curvature, and
total surface area for optimal surfaces of the Schwartz P family:
volume fraction versus mean curvature (left) and
total area per unit cell versus volume fraction (right).  
Our results (circles) are compared with those obtained by 
Anderson~\etal~(solid lines).}
\label{fig:P_volfrac_curv_area}
\end{figure}

As the Theorem~\ref{thm:extrema_total_surface_area} indicates, locally
optimal surfaces should have constant mean curvature. Because each surface
in the Schwartz P family is locally optimal, they all constant mean curvature.  
In Table~\ref{table:1}, we report the numerical values of the 
mean curvature and total surface area for the surfaces in the Schwartz P 
surface family that we computed.
\begin{table}[htb!]
\caption{The mean curvature and total surface area of the Schwartz P surface
family for different values of the volume fraction.}
\label{table:1}
\begin{tabular}{ccc} \hline
\makebox[1.20in]{Volume fraction $\volfrac_o$} & \makebox[0.95in]{$H$} & 
\makebox[0.95in]{$\cA$} \\ 
\hline
 0.25&  $ 1.67$  & 2.00 \\ 
 0.3&   $ 1.12$  & 2.14 \\ 
 0.35&  $ 0.77$  & 2.23 \\ 
 0.4&   $ 0.50$  & 2.30 \\ 
 0.45&  $ 0.24$  & 2.33 \\ 
 0.5&   $ 0.00$  & 2.34 \\ 
 0.55&  $-0.24$  & 2.33 \\ 
 0.6&   $-0.49$  & 2.30 \\ 
 0.65&  $-0.78$  & 2.23 \\ 
 0.7&   $-1.12$  & 2.14 \\ 
 0.75&  $-1.67$  & 2.00 \\ 
\hline
\end{tabular}
\end{table}
Figure~\ref{fig:P_volfrac_curv_area} compares our results with those
obtained by Anderson \etal~\footnote{The curves for 
the data obtained by Anderson \etal~are reproduced with the permission of 
the authors and were generated by digitizing the plots in~\cite{An90} using
the Plot Digitizer program.}.
Note that the data we produce matches the results of Anderson~\etal~
reasonably well.  The small discrepancy between the two sets of data is 
a consequence of (1) the lower effective grid resolution of the calculations
by Anderson~\etal~for the Schwartz P surface and (2) the digitization 
procedure used to extract the curves of Anderson~\etal~from a printed version 
of their paper. 
Figure~\ref{fig:P_volfrac_curv_area} shows that the volume is almost linearly
related to mean curvature in the volume fraction range between 
$\volfrac_o = 0.35$ and $\volfrac_o = 0.65$. Beyond this range, the 
mean curvature varies more rapidly as a function of the volume fraction.
Anderson's data shows that near volume fractions of 0.75 and 0.25, the 
relationship between the volume fraction and mean curvature reverses 
direction.  Due to the limitations of our computational method, we have not
yet been able to reproduce the portions of curve between points B and C.
We conjecture that our numerical method fails for these portions of the 
curves because they are not local minima of the surface area.

\subsubsection{Schwartz D surface family}
\label{subsubsection5.3.2}
The Schwartz D surface has largest surface area among the three triply 
periodic surfaces studied in this paper.  Since large surface area is 
inversely correlated with minimum distance between distinct sheets of
a surface, the calculations for the Schwartz D surface required 
higher grid resolution than those for the Schwartz P surface.  For
the present calculations, we employed a grid of size 
$250 \times 250 \times 250$. Coarser meshes (\eg $150 \times 150 \times 150$) 
failed at extreme values of the volume fractions 
(\eg $\volfrac_o = 0.15, 0.85$) because the optimal surfaces at these volume 
fractions have very narrow connections.   Due to the larger computational
grid, the run time for these simulations took significantly longer than
for the Schwartz P-type surfaces.
Reaching an area stopping criterion of $\Delta \cA_n < 10^{-5}$ took 
approximately 30 hours for a 4 node parallel calculation on the Opteron cluster cluster.  As with the Schwartz P-type simulations, the running time depended
on the volume fraction constraints.
\begin{figure}[htb!]
\centerline{
\includegraphics[width=1.6in]{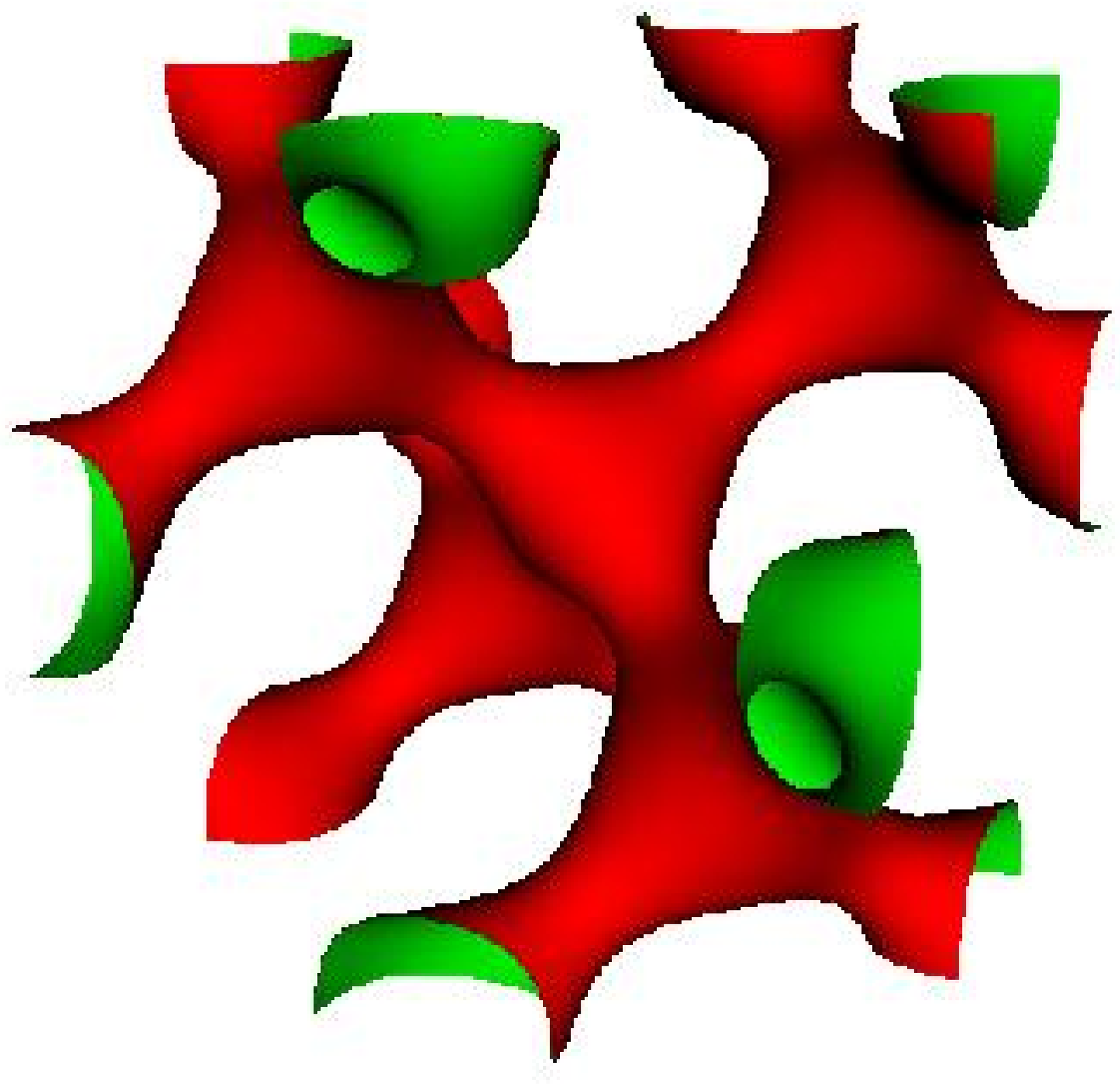} \hspace{0.0in}
\includegraphics[width=1.6in]{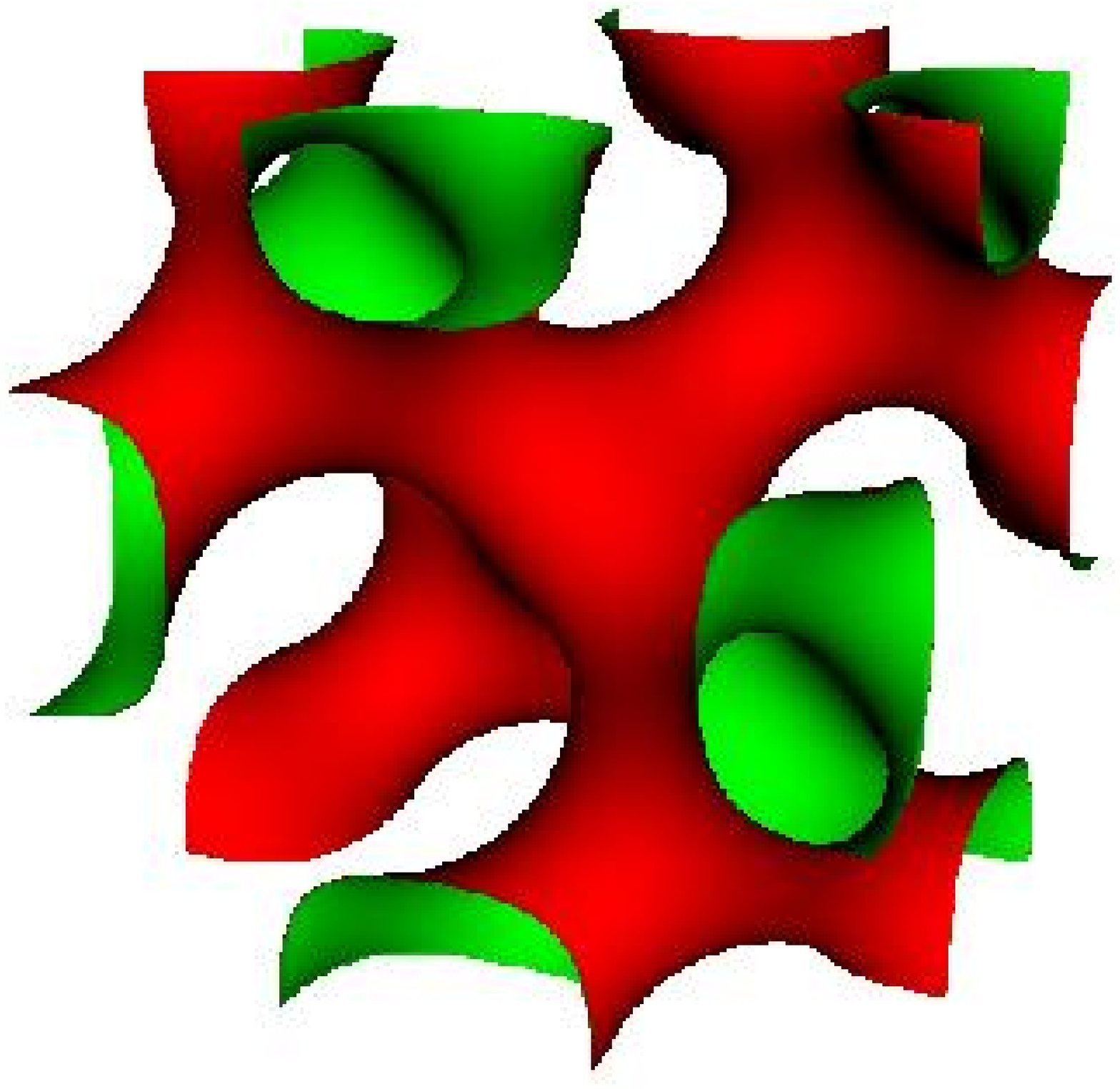} \hspace{0.0in}
\includegraphics[width=1.6in]{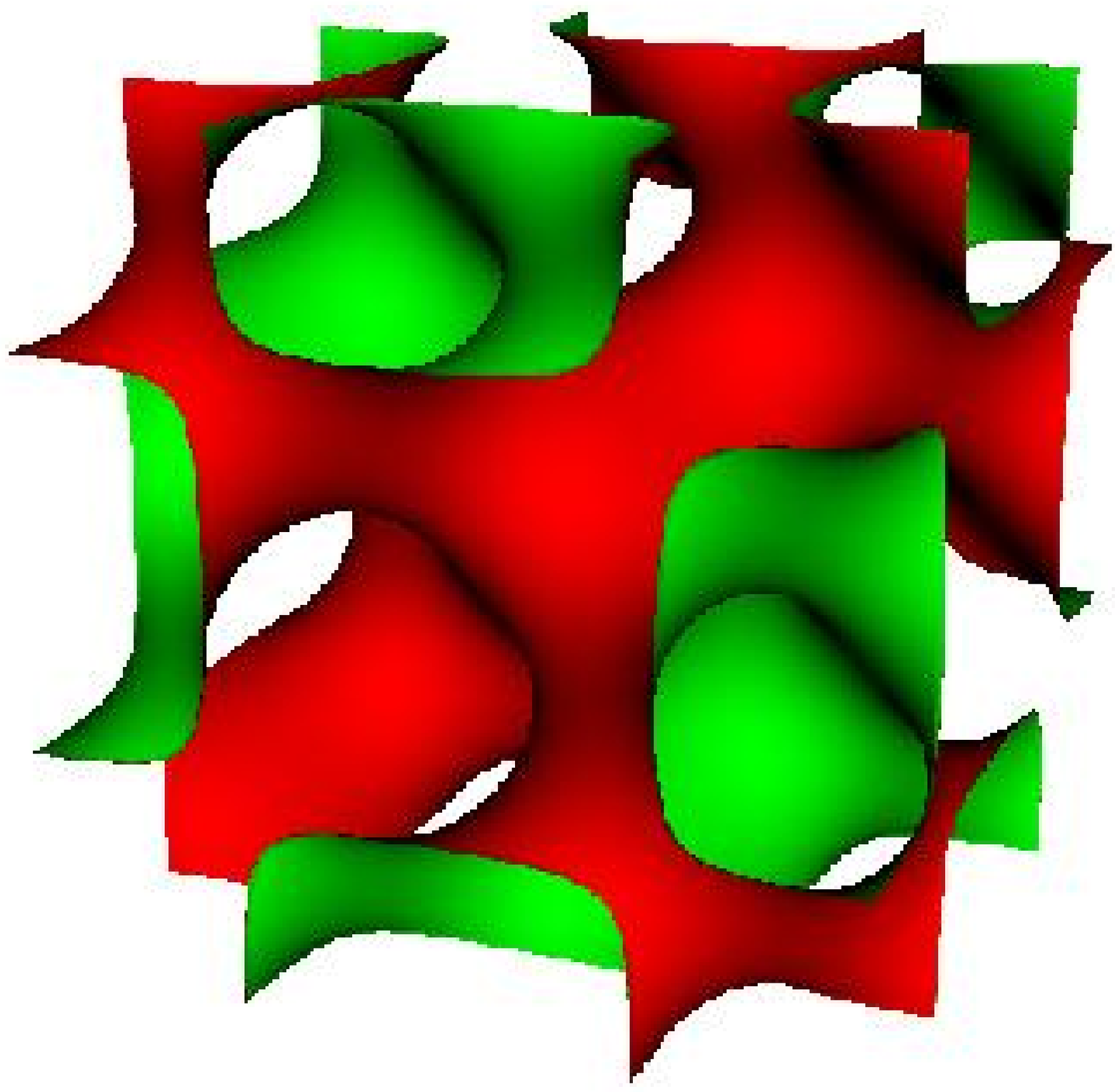}
}
\vspace{0.1in}
\centerline{
\includegraphics[width=1.6in]{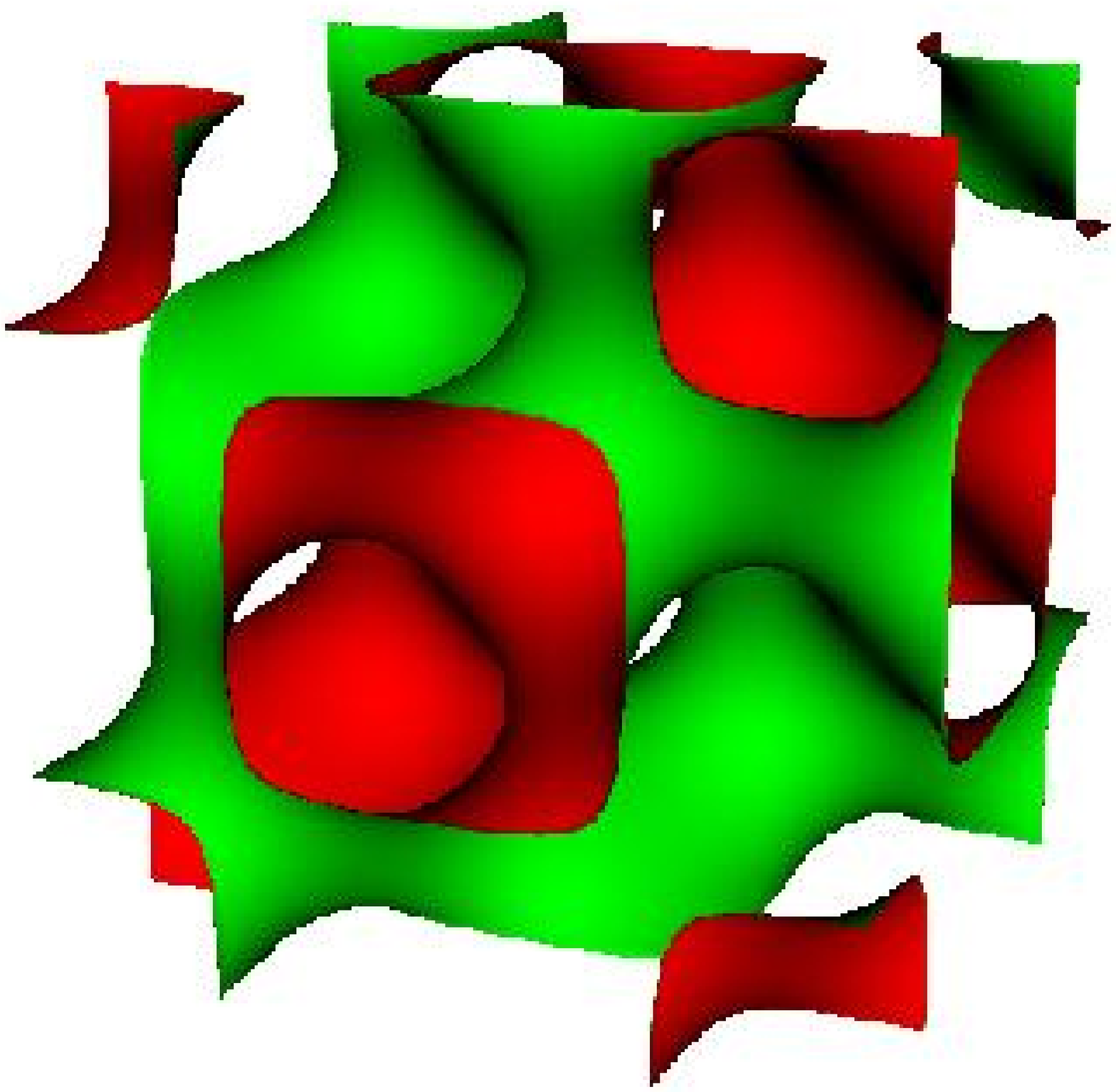} \hspace{0.0in}
\includegraphics[width=1.6in]{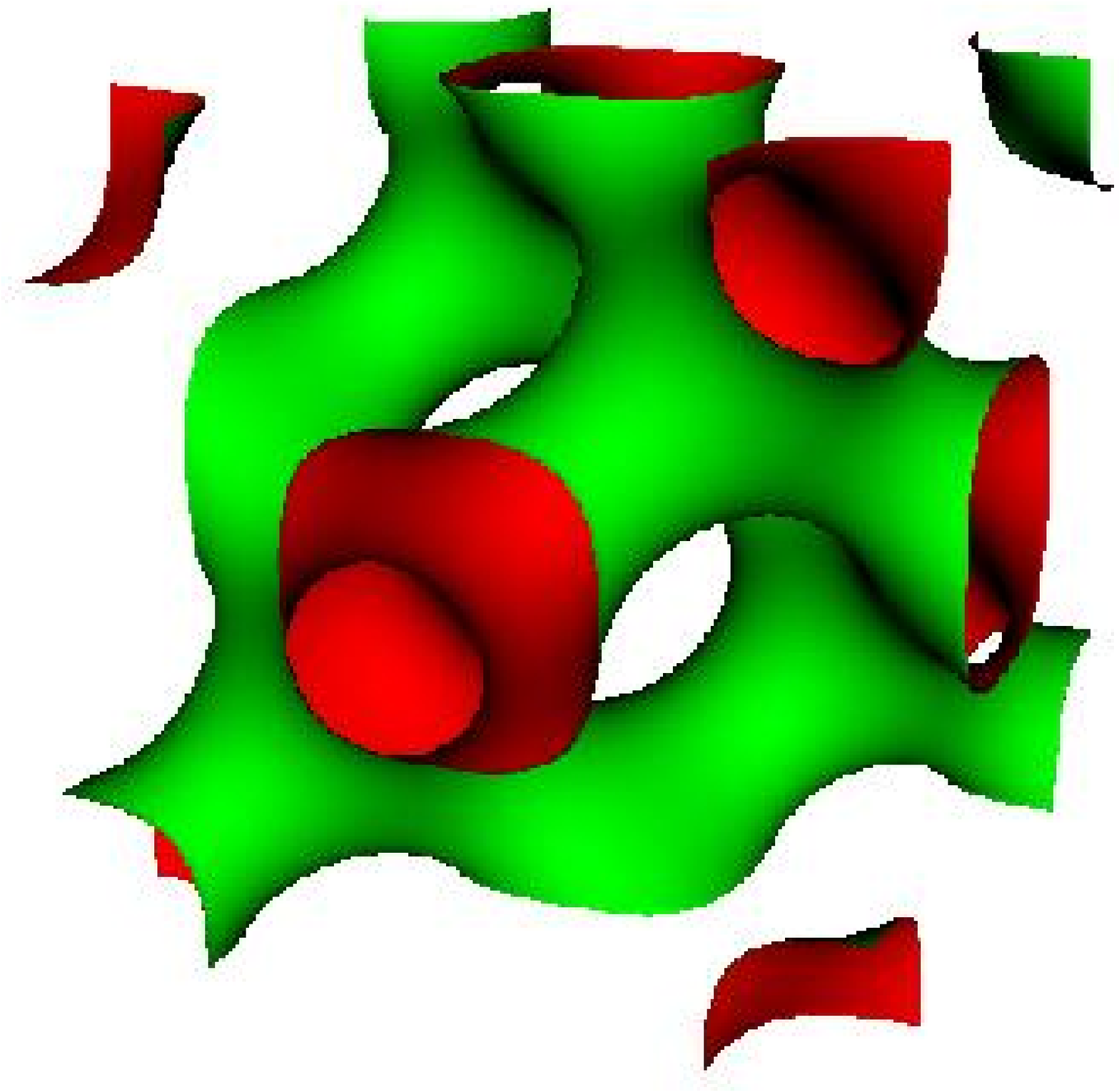} \hspace{0.0in}
\includegraphics[width=1.6in]{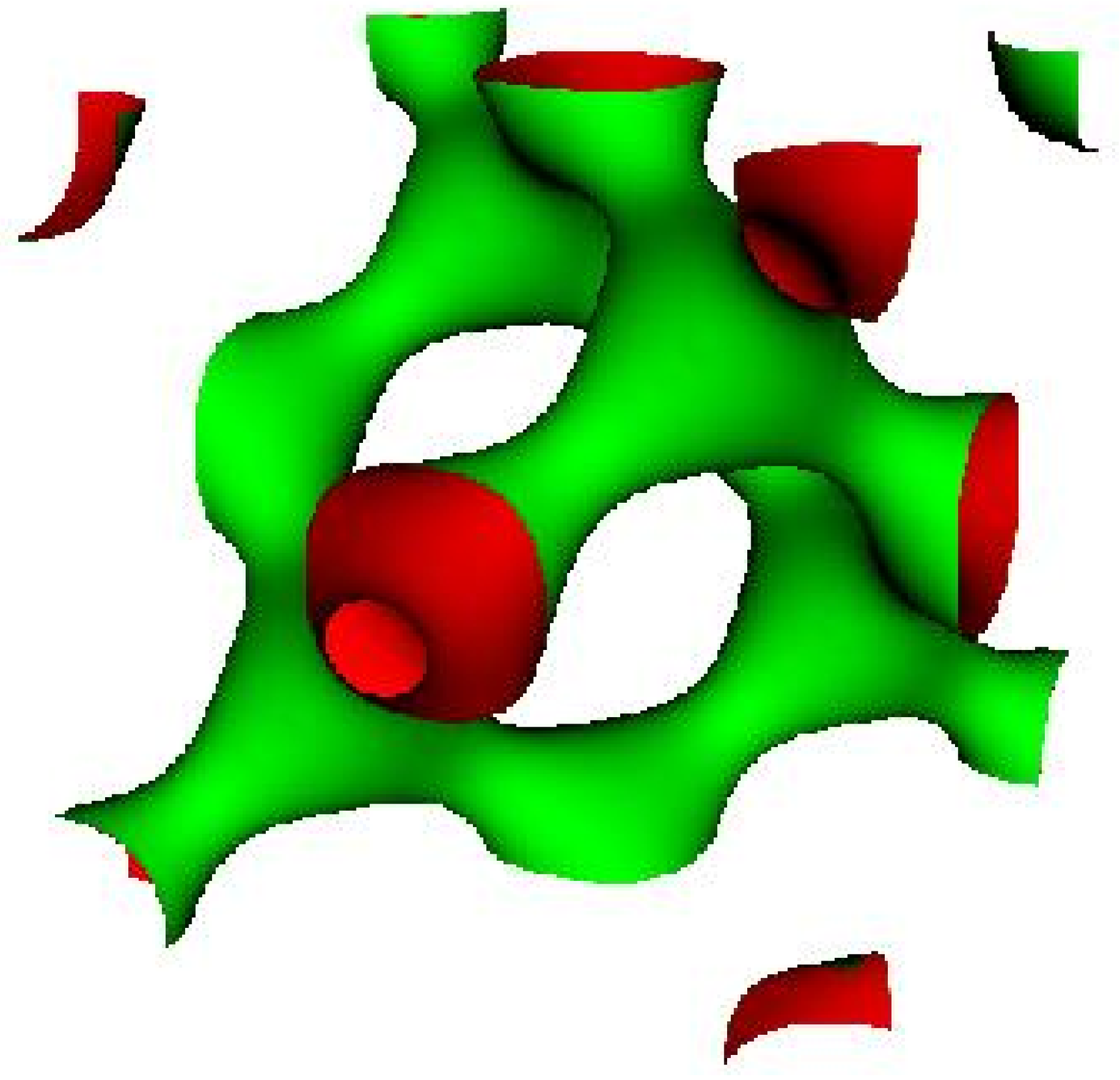}
}
\caption{Family of surfaces possessing the symmetry and 
topology of the Schwartz D surface but with different volume fractions for 
phase 1. 
Upper left to upper right panel: $\volfrac_o = 0.15, 0.3, 0.45$. 
Lower left to lower right panel: $\volfrac_o = 0.55, 0.7, 0.85$.}
\label{fig:D_different_volume}
\end{figure}

\begin{figure}[htb!]
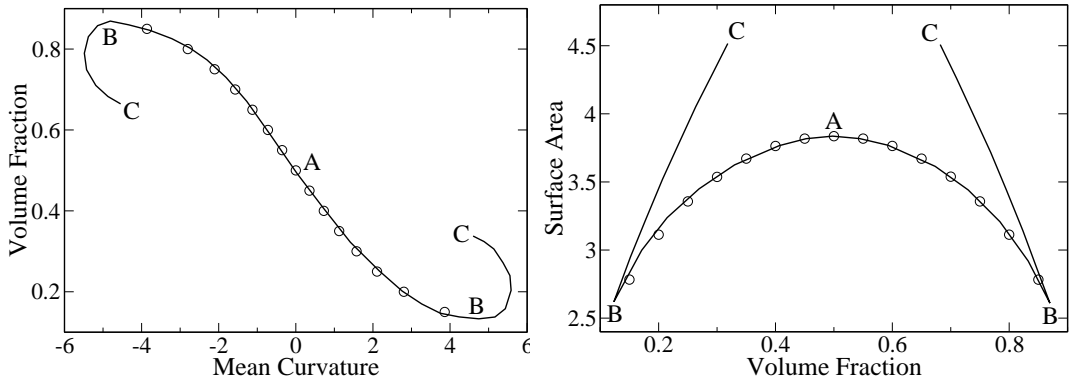

\vspace{0.1in}
\centerline{
\includegraphics[width=2.75in]{figs/D_curv_vol_250}
\includegraphics[width=2.75in]{figs/D_vol_area_250}
}
\caption{Relationships between the volume fraction, mean curvature, and
total surface area for optimal surfaces of the Schwartz D family:
volume fraction versus mean curvature (left) and
total area per unit cell versus volume fraction (right).
Our results (circles) are compared with those obtained by 
Anderson~\etal~(solid lines).}
\label{fig:D_volfrac_curv_area}
\end{figure}

Figure~\ref{fig:D_volfrac_curv_area} shows excellent agreement with 
the results of Anderson~\etal~\footnote{In the calculations of Anderson~\etal,
the effective grid resolution for the Schwartz D surface was higher than
for the Schwartz P surface.}.
However, as with the Schwartz P-type surfaces,
our algorithm cannot compute surfaces between points B and C on the
curve.  The numerical values of the mean curvature and total surface area for
the surfaces in the Schwartz D surface family are reported in 
Table~\ref{table:2}.

\begin{table}[htb!]
\caption{The mean curvature and total surface area of the Schwartz D surface
family for different values of the volume fraction. }
\label{table:2}
\begin{tabular}{ccc} \hline
\makebox[1.20in]{Volume fraction $\volfrac_o$} & \makebox[0.95in]{$H$} & 
\makebox[0.95in]{$\cA$} \\ 
\hline
 0.15&  $ 3.86$  & 2.78 \\ 
 0.2&   $ 2.80$  & 3.11 \\ 
 0.25&  $ 2.11$  & 3.36 \\ 
 0.3&   $ 1.58$  & 3.54 \\ 
 0.35&  $ 1.12$  & 3.67 \\ 
 0.4&   $ 0.73$  & 3.76 \\ 
 0.45&  $ 0.35$  & 3.82 \\ 
 0.5&   $ 0.00$  & 3.84 \\ 
 0.55&  $-0.35$  & 3.82 \\ 
 0.6&   $-0.72$  & 3.76 \\ 
 0.65&  $-1.13$  & 3.67 \\ 
 0.7&   $-1.57$  & 3.54 \\ 
 0.75&  $-2.10$  & 3.36 \\ 
 0.8&   $-2.80$  & 3.11 \\ 
 0.85&  $-3.86$  & 2.78 \\ 
\hline
\end{tabular}
\end{table}

\subsubsection{Schoen G surface family}
\label{subsubsection5.3.3}
\begin{figure}[htb!]
\centerline{
\includegraphics[width=1.6in]{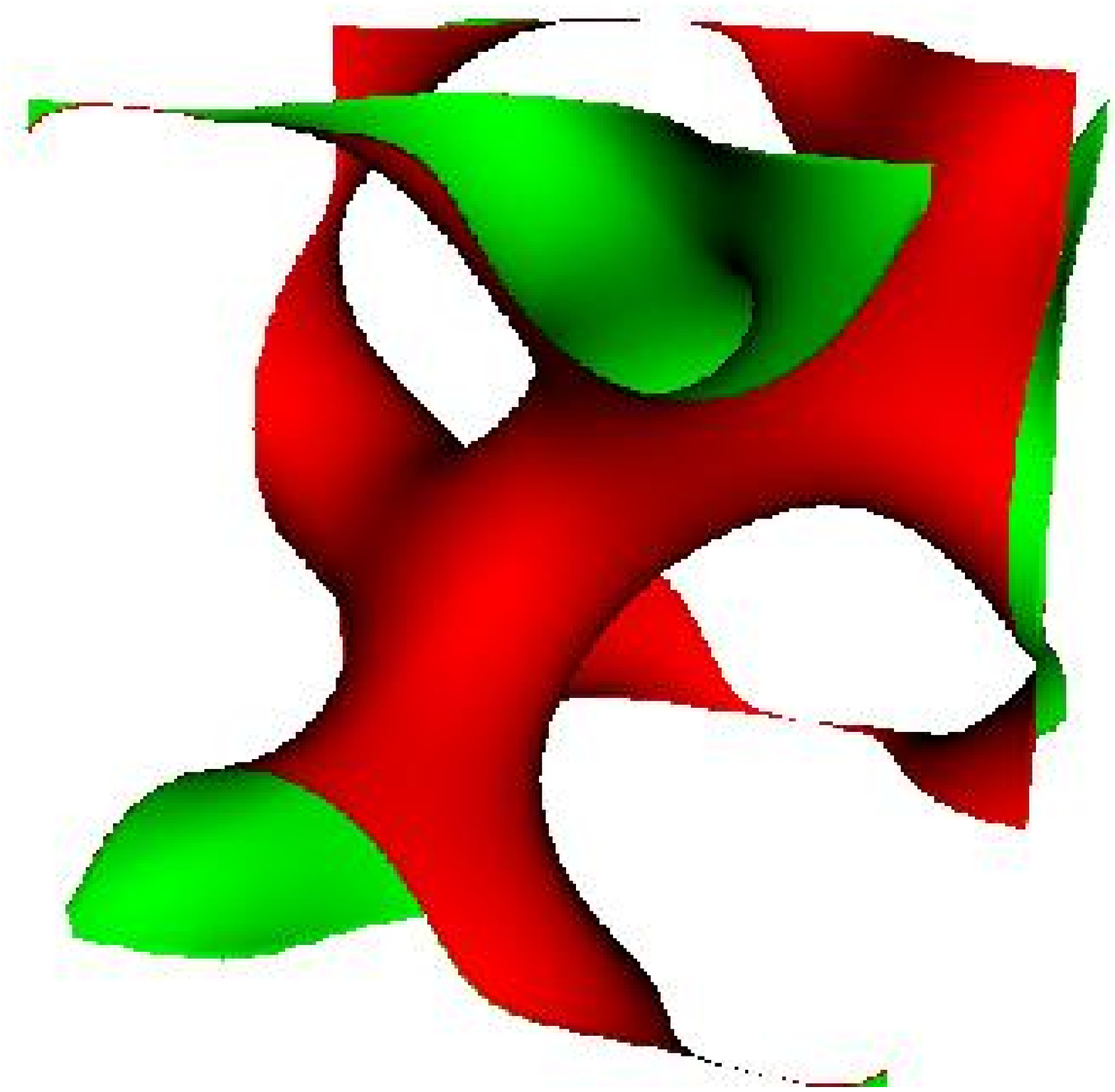} \hspace{0.0in}
\includegraphics[width=1.6in]{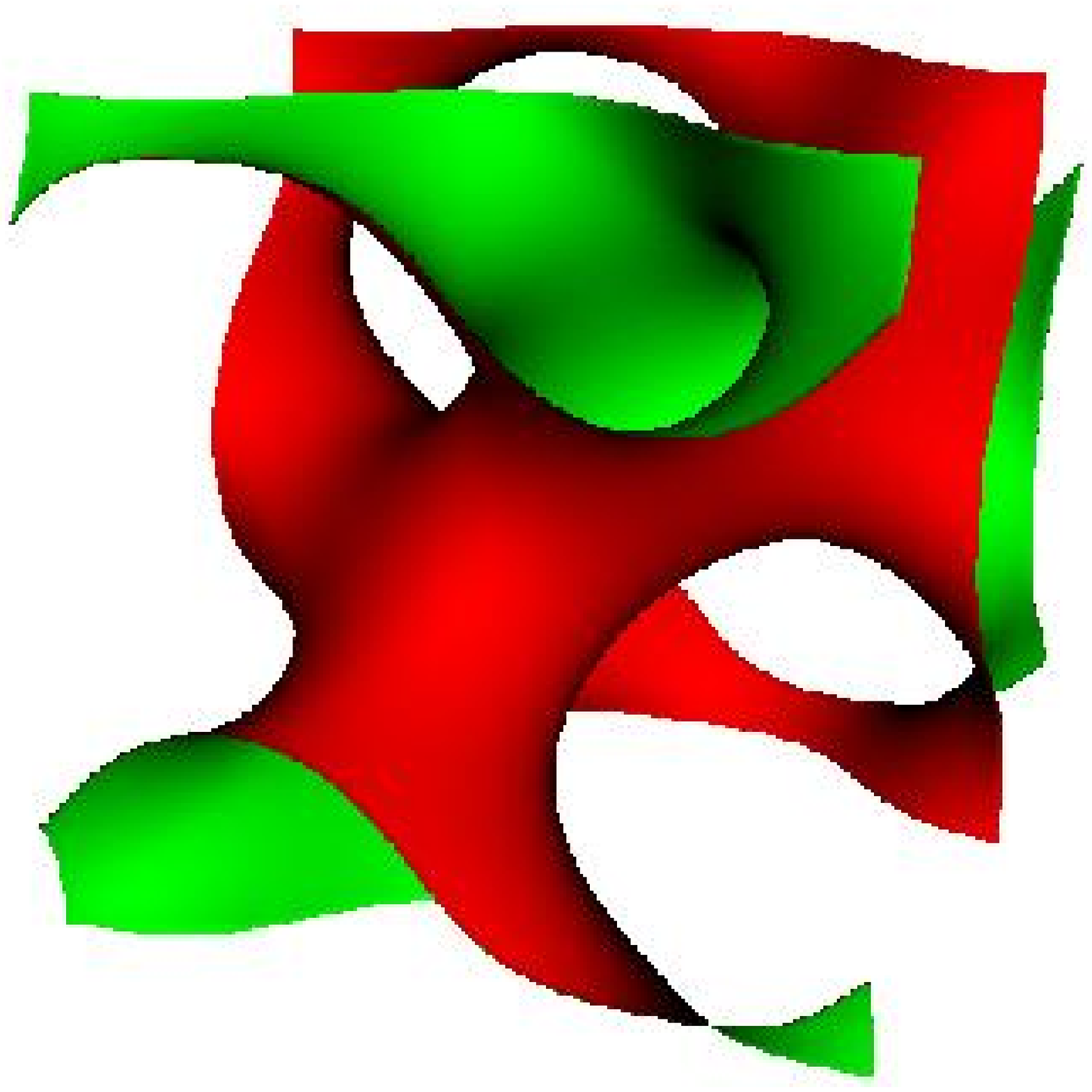} \hspace{0.0in}
\includegraphics[width=1.6in]{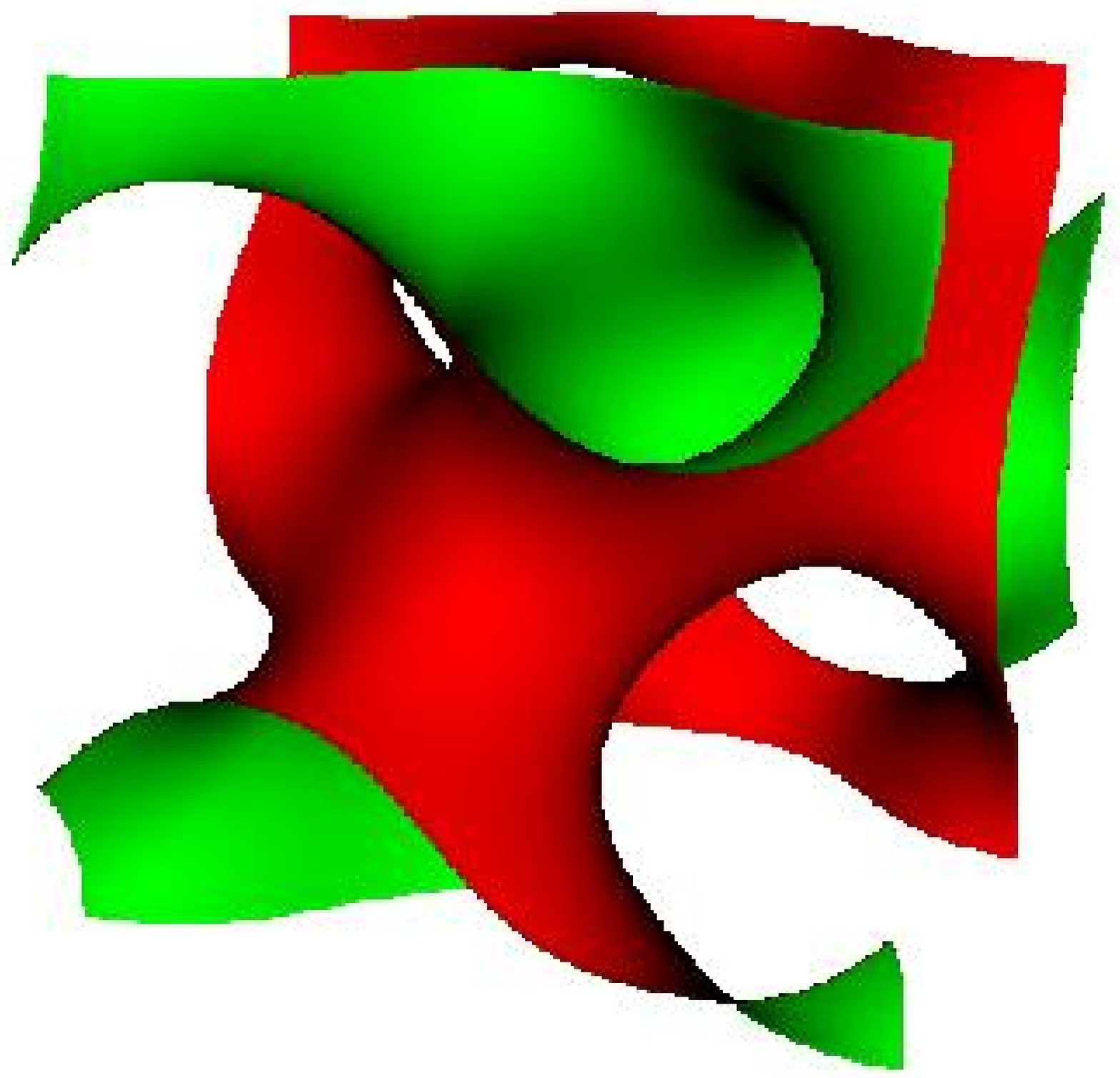}
}
\vspace{0.3in}
\centerline{
\includegraphics[width=1.6in]{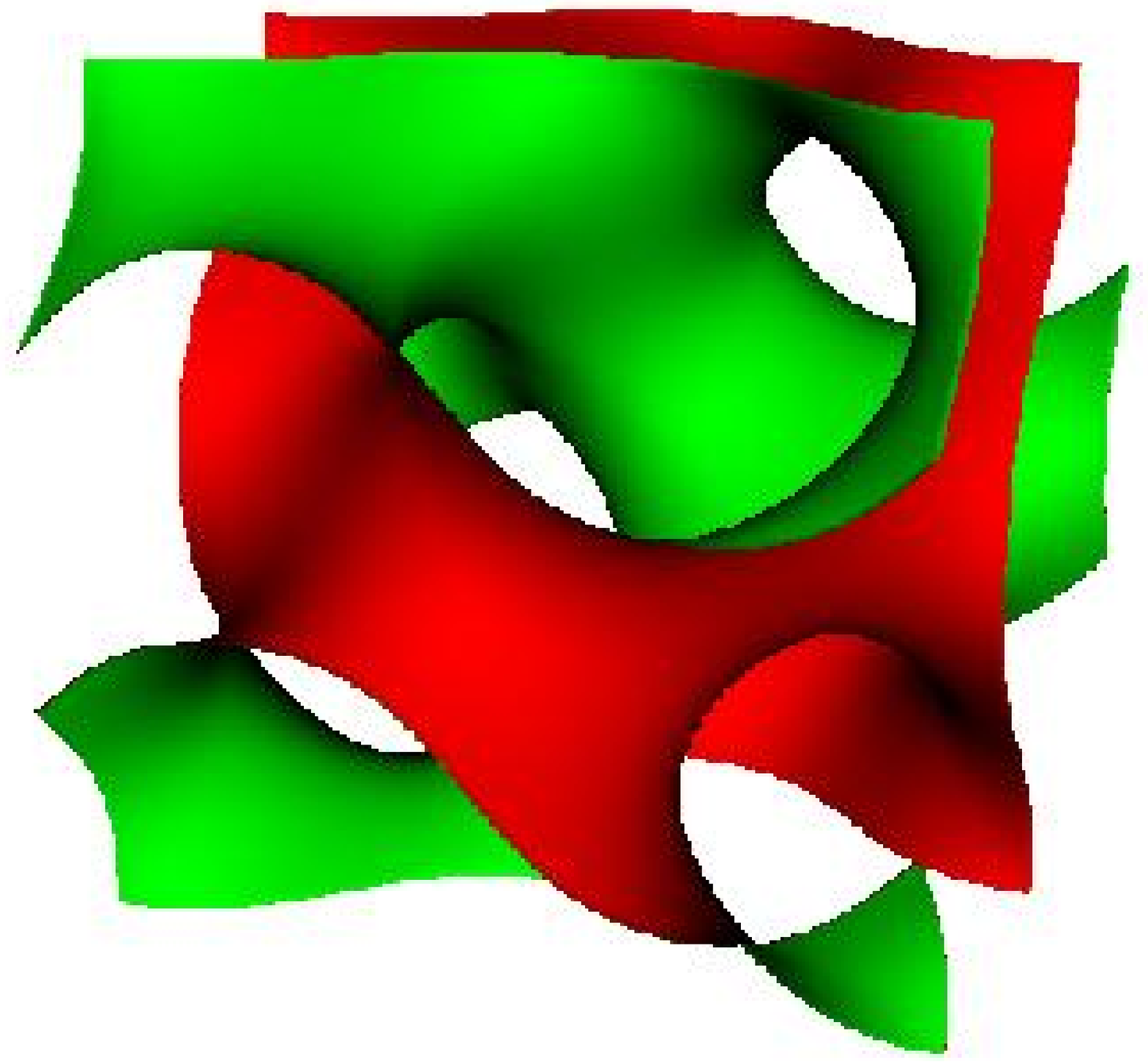} \hspace{0.0in}
\includegraphics[width=1.6in]{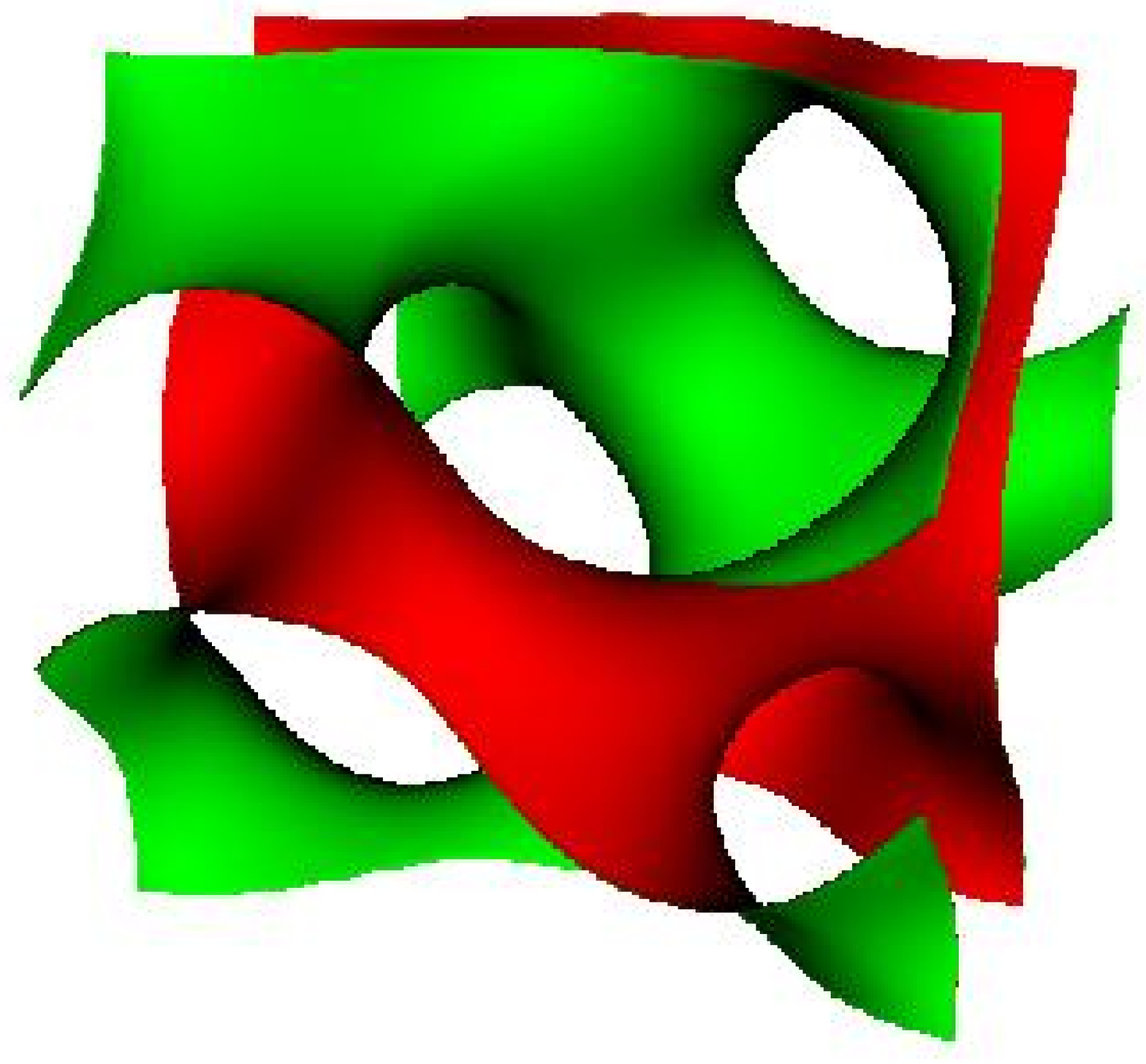} \hspace{0.0in}
\includegraphics[width=1.6in]{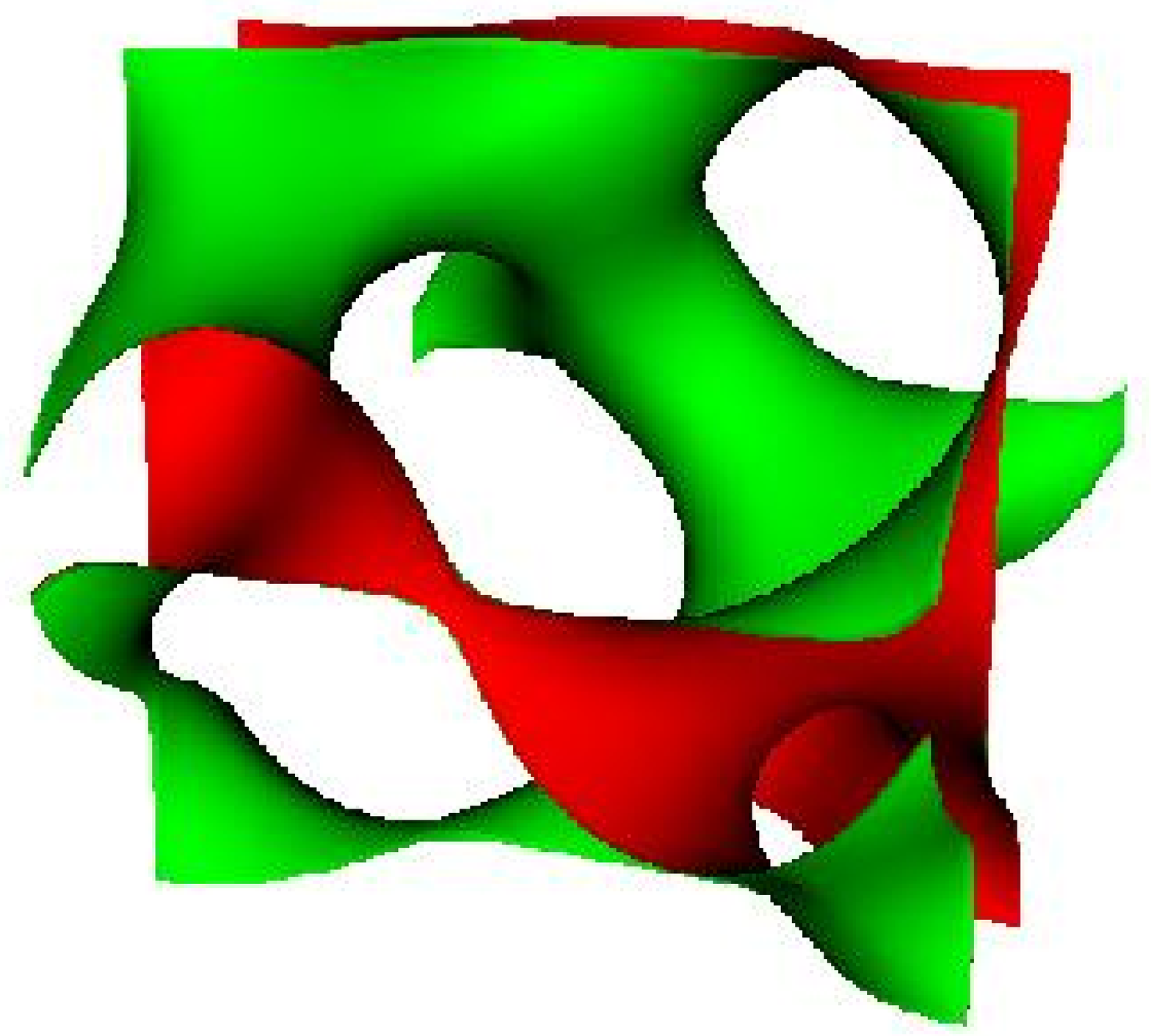}
}
\caption{Family of surfaces possessing the symmetry and topology
of the Schoen G surface but with different volume fractions for phase 1. 
Upper left to upper right panel: $\volfrac_o = 0.2, 0.3, 0.4$. 
Lower left to lower right panel: $\volfrac_o = 0.6, 0.7, 0.8$.}
\label{fig:G_different_volume}
\end{figure}

\begin{figure}[htb!]
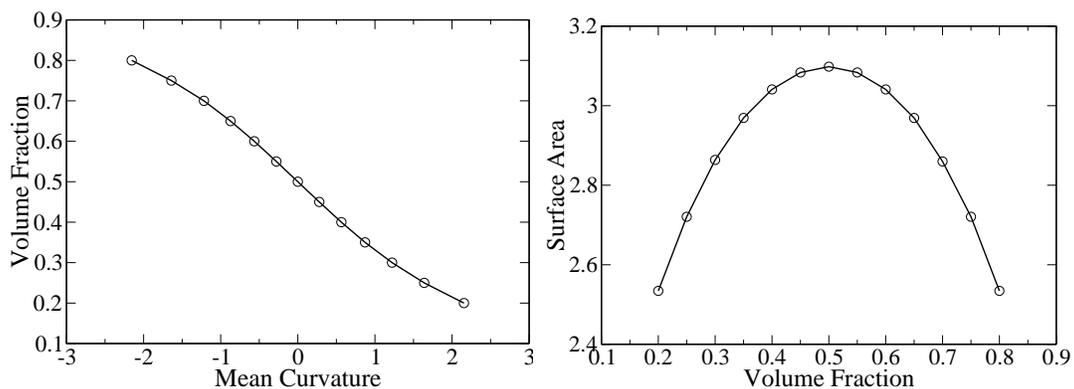

\centerline{
\includegraphics[width=2.75in]{figs/G_curv_vol_200}
\includegraphics[width=2.75in]{figs/G_vol_area_200}
}
\caption{Relationships between the volume fraction, mean curvature, and
total surface area for optimal surfaces of the Schoen G family:
volume fraction versus mean curvature (left) and
total area per unit cell versus volume fraction (right).}
\label{fig:G_volfrac_curv_area}
\end{figure}

For these calculations, we employed a grid of size $200 \times 200 \times 200$. 
Because Anderson \etal~did not study Schoen G surfaces, we explored the valid 
range of volume fractions by stepping the volume fraction in increments
of 0.05 from a volume fraction of $0.5$ until we observed optimal surfaces
that no longer possessed the symmetry of the Schoen G family 
(Figure~\ref{fig:G_different_volume}).
We detected the symmetry change (as well as sudden changes in the 
mean curvature and total surface area) beyond $\volfrac_o = 0.8$ and 
$\volfrac = 0.2$.  These results were consistent with simulations on finer 
meshes (\eg $250 \times 250 \times 250$), which leads us to conclude that the 
$\volfrac_o = 0.8$ and $\volfrac_o = 0.2$ are close to the turning points 
in the volume fraction versus mean curvature graph 
(Figure~\ref{fig:G_volfrac_curv_area}).

As shown in Figure~\ref{fig:G_volfrac_curv_area}, the mean curvature range 
for the Schoen G family of surfaces is similar to that of the Schwartz P 
family of surfaces in Section~\ref{subsubsection5.3.1}. However, the total 
surface area of Schoen G-type surfaces are bigger than those of the 
Schwartz P-type surfaces with the the same volume fractions.
Numerical values for the mean curvature and total surface area for 
surfaces in the Schoen G family are reported in Table~\ref{table:3}.

\begin{table}[htb!]
\caption{The mean curvature and total surface area of the Schoen G surface 
family for different values of the volume fraction. }
\label{table:3}
\begin{tabular}{ccc} \hline
\makebox[1.20in]{Volume fraction $\volfrac_o$} & \makebox[0.95in]{$H$} & 
\makebox[0.95in]{$\cA$} \\ 
\hline
 0.2&   $ 2.15$  & 2.53 \\ 
 0.25&  $ 1.64$  & 2.72 \\ 
 0.3&   $ 1.22$  & 2.86 \\ 
 0.35&  $ 0.87$  & 2.97 \\ 
 0.4&   $ 0.56$  & 3.04 \\ 
 0.45&  $ 0.28$  & 3.08 \\ 
 0.5&   $ 0.00$  & 3.10 \\ 
 0.55&  $-0.28$  & 3.08 \\ 
 0.6&   $-0.56$  & 3.04 \\ 
 0.65&  $-0.87$  & 2.97 \\ 
 0.7&   $-1.22$  & 2.86 \\ 
 0.75&  $-1.64$  & 2.72 \\ 
 0.8&   $-2.15$  & 2.53 \\ 
\hline
\end{tabular}
\end{table}

Following Anderson \etal~, it would be interesting to complete the graphs of 
the volume fraction versus mean curvature and total surface area. 
However, as with the Schwartz P and D surfaces, we are currently unable to
generate surfaces beyond the turning points in the graph of the volume 
fraction versus mean curvature.

\section{Conclusions}
\label{section6}
In this paper, we have generated and probed the structure of the constant
mean curvature, triply periodic surfaces (including minimal surfaces) that 
arise when the total surface area is minimized subject to a constraint on the 
volume fraction of the regions that the surface separates.  
Specifically, we have shown that the Schwartz P, Schwartz D, and Schoen G
minimal surfaces are local minima of the total surface area under a volume 
fraction constraint, studied the properties of families of surfaces 
possessing the symmetry of the Schwartz P, Schwartz D, and Schoen G surfaces 
at varying volume fractions of the phases separated by the surface, and
generated new minimal surfaces (e.g. with high genus) as well
as new surfaces with constant, nonzero mean curvature (e.g. family of
surfaces with the symmetry of the Schoen G surface).  
Unlike many studies on surfaces, our perspective is interesting 
because it draws attention to the properties of surfaces when global 
geometric constraints rather than local differential conditions are imposed.  
Many of our results complement and extend the work by Andersen 
\etal~\cite{An90} on nonzero, constant mean curvature surfaces. 
We emphasize that an important feature of our approach is our ability 
to easily control the volume fraction of the phases.

The basic approach that we have taken to study optimal triply periodic
surfaces has been to use ideas from the level set method and shape 
optimization.  This well-known framework has allowed us to formulate the 
problem in a manner that is amenable both to theoretical analysis and numerical 
simulation.  Theoretically, it allowed us to prove that surfaces with 
constant mean curvature exactly are precisely the ones that optimize 
total surface area under a volume fraction constraint on the phases.  
Computationally, the framework naturally led (via a steepest descent 
procedure) to a generalization of the two-dimensional shape optimization 
algorithms presented in~\cite{osher01} and \cite{alexandrov_2005}. 
We note that our discussion of the optimization algorithm provides a
more comprehensive description of how variational calculus is used in
the context of shape optimization and level set methods.
In particular, we demonstrate how the variational formulation of the 
problem can be used to draw theoretical conclusions (in addition to 
forming the foundation for the computational method).

Using our numerical algorithm, we explored several optimal triply periodic 
surfaces.  When the volume fractions of the two phases are equal, our
theorem implies that the Schwartz P, Schwartz D, and Schoen G 
surfaces are all local optima of the surface area;  numerical 
simulation shows that these surfaces are actually local minima of the
total surface area.  For unequal volume fractions of the phases, we 
partially reproduced the volume fraction versus mean curvature and 
total surface area versus volume fraction results obtained by 
Anderson~\etal~\cite{An90} for the Schwartz P and D surfaces.  
Unfortunately, due to limitations in our algorithm, we were unable to
reproduce the results for extreme values of the mean curvature.  We 
conjecture the surfaces at extreme values of the mean curvature 
may no longer be local minima of the total surface area.  This question
certainly deserves require further investigation.
In addition to studying the Schwartz P and Schwartz D families of surfaces, 
we also extended the work of Anderson \etal~by determining the relationship 
between the volume fraction, mean curvature, and surface area for the Schoen 
G surface.

One of the original goals our investigation was to show that the triply
periodic surface with minimal total surface area at a volume fraction of 
1/2 is the Schwartz P surface.  While our numerical simulations provided
us with empirical evidence supporting this conjecture, we could not prove
these results because our numerical algorithm and theoretical analysis are 
limited to information about local minima of the total surface area.  
Thus, the question of the \emph{global} optimality of the Schwartz P surface
remains an open question for future investigation.

\subsection{Future directions}
To close, we mention a few future directions for our work.  First,
there are research questions that might allow us to gain more insight 
into the question of global minimality of the Schwartz P surface.  
For instance, a way to quantify the width of the basin of attractions for the 
local minima of the total surface area could allow us to be more systematic in 
our search through the space of triply periodic surfaces (which could be
carried out using our numerical optimization algorithm.) 
Our numerical algorithm could also be used to explore the optimal 
surfaces possessing a wider array of symmetries.

Second, with a working implementation of our numerical algorithm, we are now
in a position to study important questions about the physical properties of 
two-phase materials whose phases are separated by optimal surfaces.  
For instance, it has been conjectured that microstructures where two phases 
of unequal volume fraction are separated by a surface that minimize the 
total surface area lie on the upper bound of effective transport properties
for composite materials~\cite{To02d,To03a}. One obstacle to collecting 
evidence supporting or refuting this conjecture has been our ability to 
generate surfaces of minimal surface area.  Our optimization algorithm
effectively removes this barrier and makes the multifunctional optimization
problem for unequal volume fractions more tractable.

\section*{Acknowledgments}
YJ and ST gratefully acknowledge the support of the Air Force Office for 
Scientific Research under grant number F49620-03-1-0406.
KTC gratefully acknowledges the support of the National Science Foundation
through grant DMR-0502946.
The authors thank H.~T.~Davis for his help in our attempts to obtain
the original data for the plots from~\cite{An90}, D.~J.~Srolovitz for
his encouragement and support in the development of the LSMLIB parallel 
level set method software library, A. Donev for helpful comments on the
manuscript, and Bill Wischer for his technical support in carrying out 
the parallel computations.


\end{document}